\documentclass[11pt,amsfonts,fullpage]{amsart}

\usepackage{fullpage}

\newtheorem{theorem}{Theorem}[section]

\newtheorem{corollary}[theorem]{Corollary}

\newtheorem{definition}[theorem]{Definition}

\newtheorem{lemma}[theorem]{Lemma}

\newtheorem{proposition}[theorem]{Proposition}

\newtheorem{remark}[theorem]{Remark}

\newtheorem{conjecture}[theorem]{Conjecture}

\newtheorem{example}[theorem]{Example}

\def\endproof{\qed \medskip}
\def\blacksquare{\hbox to .60em{\vrule width .60em height .60em}}

\begin{document}

\title[ ]{Canonical Metrics on 3-Manifolds and 4-Manifolds}

\author[ ]{Michael T. Anderson}

\thanks{Partially supported by NSF Grant DMS 0305865}

\maketitle

\begin{flushright}
In Memory of S.S. Chern
\end{flushright}

\setcounter{section}{0}

\section{Introduction.}
\setcounter{equation}{0}

 In this paper, we discuss recent progress on the existence of canonical 
metrics on manifolds in dimensions 3 and 4, and the structure of moduli spaces 
of such metrics. The existence of a ``best possible'' metric on a given closed 
manifold is a classical question in Riemannian geometry, attributed variously to 
H. Hopf and R. Thom, see [22] for an interesting perspective. A good deal of motivation 
for this question comes from the case of surfaces; the uniformization theorem in 
dimension 2 has a multitude of consequences in mathematics and physics. Further, 
there are strong reasons showing that the closest relations between geometry and 
topology occur in dimensions 2, 3 and 4.

 The precise formulation of the question in dimension 3 is given by
Thurston's Geometrization Conjecture. This conjecture describes completely when 
a given 3-manifold admits a canonical metric (defined to be a metric of constant 
curvature or more generally a locally homogeneous metric), and thus determines 
exactly what the obstructions are to the existence of such a metric. Moreover, 
it describes how an arbitrary 3-manifold decomposes into topologically 
essential pieces, each of which admits a canonical metric, resulting in the 
topological classification of 3-manifolds. The apparent solution of the 
Geometrization Conjecture by Perelman is one of the most spectacular breakthroughs 
in geometry and topology in the past several decades.

 The Thurston picture will be reviewed in more detail in \S 2, for the light it 
sheds on what might be hoped for or expected in dimension 4. Since 
there is already considerable analysis and discussion of the details of 
Perelman's work elsewhere, we will not discuss this in any detail 
here. We do however give one application of his work, (since this does 
not seem to be widely known), namely the determination of the value of the 
Yamabe invariant or Sigma constant $\sigma (M)$ of all 3-manifolds $M$ for 
which $\sigma (M) \leq 0$, cf. \S 2.

\medskip

 Thus, the bulk of the paper concerns dimension 4. Canonical metrics 
will be defined to be metrics minimizing, (or possibly just critical 
points for), one of the classical and natural curvature functionals 
${\mathcal F}$ on the space of metrics ${\mathbb M}$ on a given oriented 4-manifold $M$:
\begin{equation} \label{e1.1}
{\mathcal R}^{2}, \ {\mathcal W}^{2}, \ {\mathcal W}_{+}^{2}, \ {\mathcal W}_{-}^{2}, 
\ {\mathcal R}ic^{2}. 
\end{equation}
These are respectively the square of $L^{2}$ norm of the Riemann 
curvature $R$, Weyl curvature $W$, its self-dual and anti-self-dual 
components, $W_{+}$, $W_{-}$, and Ricci curvature $Ric$. We will also 
consider, but in much less detail, the scalar curvature functionals
\begin{equation} \label{e1.2}
{\mathcal S}^{2}, -{\mathcal S}|_{{\mathcal Y}}, 
\end{equation}
given by the square of the $L^{2}$ norm of the scalar curvature $s$, and 
the restriction of the total scalar curvature to the space ${\mathcal Y}$ 
of unit volume Yamabe metrics on $M$. The Chern-Gauss-Bonnet theorem [17] 
relating the functionals in (1.1)-(1.2) plays a crucial r\^ole in the analysis to 
follow. 

 Einstein metrics, satisfying 
\begin{equation} \label{e1.3}
Ric_{g} = \frac{s}{n}g, 
\end{equation}
$n = dim M$, are critical points of all of the functionals in 
(1.1)-(1.2), and in many cases, Einstein metrics are minimizers. 
However, there are large classes of minimizers of ${\mathcal W}^{2}$, or 
the related ${\mathcal W}_{+}^{2}$, ${\mathcal W}_{-}^{2}$ which are not 
Einstein. Critical metrics of ${\mathcal W}^{2}$ are Bach-flat metrics, 
satisfying the Bach equations, and include conformally flat as well as 
half-conformally flat (self-dual or anti-self-dual) metrics and these 
classes of metrics are also often minimizers. It does not seem to be known 
if there are any other minimizers, or even critical points of the functionals 
in (1.1) or (1.2), which are not Einstein or half-conformally flat. 

 Just as in dimension 3, on a general 4-manifold, metrics minimizing a 
particular functional in (1.1) will not exist. Until relatively recently, 
the only known obstructions to the existence of Einstein metrics were the 
Hitchin-Thorpe inequality $\chi(M) \geq \frac{3}{2}|\tau(M)|$ between the Euler 
characteristic and signature of $M$, and Gromov's improvement of this, based on 
the simplicial volume. In the past decade or so, many further obstructions have 
been found by LeBrun and others, which show for instance that the existence of 
an Einstein metric on a 4-manifold $M$ often depends strongly on the smooth 
structure of $M$ as opposed to just the topological structure; we refer to 
[28], [29] for overviews of the current status of these issues. Nevertheless, 
one is far from having a comprehensive understanding of the obstructions to 
Einstein or half-conformally compact metrics on a given 4-manifold. 

\medskip

 In \S 3, we survey in some detail the currently known results 
concerning the structure of the moduli space of Einstein metrics and 
moduli spaces ${\mathcal M}_{{\mathcal F}}$ of the functionals ${\mathcal F}$ in 
(1.1). 

 In \S 4, these results, and the methods used in their proof, are extended 
to prove a general result on the weak or idealized existence of minimizers 
of the functionals in (1.1). The main result is summarized as follows, but 
we refer to Theorem 4.10, both for the definitions involved  and for a 
more precise formulation. 

\begin{theorem} \label{t 1.1.}
 Let $M$ be a closed, oriented 4-manifold and let ${\mathcal F}$ be either of 
the functionals ${\mathcal R}^{2}$ or ${\mathcal R}ic^{2}$ in (1.1). Then there 
exist minimizing sequences $\{g_{i}\}$ for ${\mathcal F}$ on the space 
${\mathbb M}_{1}$ of unit volume metrics metrics on $M$ which exhibit one 
of the following behaviors: 

 {\rm (I)}. The sequence $\{g_{i}\}$ converges in the Gromov-Hausdorff topology 
to a compact, oriented orbifold $(V, g_{0})$ associated to $M$, possibly 
reducible, with $C^{\infty}$ metric $g_{0}$ on the regular set $V_{0}$, and having 
a $C^{0}$ extension across the singular points. One has 
\begin{equation} \label{e1.4}
{\mathcal F} (g_{0}) \leq \inf_{g\in{\mathbb M}_{1}}{\mathcal F} (g), \ \ {\rm and} \ \ 
vol_{g_{0}}V = 1. 
\end{equation}

  {\rm (II)}. The sequence $\{g_{i}\}$ collapses everywhere, i.e. $inj_{g_{i}}(x) 
\rightarrow 0$, for all $x \in M$, and on the complement of a finite collection $B_{i}$ 
of arbitrarily small balls $B_{i} = \cup_{k}B_{z_{k}}(\varepsilon_{i})$, $\varepsilon_{i} 
\rightarrow 0$, the sequence $\{g_{i}\}$ collapses with locally bounded curvature along 
a sequence of F-structures on $M \setminus B_{i}$. 

 {\rm (III)}. The sequence $\{g_{i}\}$ converges in the pointed Gromov-Hausdorff topology 
to a maximal open orbifold $\Omega$, possibly reducible and possibly 
empty, with $C^{\infty}$ smooth metric $g_{0}$ on the regular set $\Omega_{0}$, 
$C^{0}$ across singular points, and satisfying
\begin{equation} \label{e1.5}
{\mathcal F} (g_{0}) \leq  \inf_{g\in{\mathbb M}_{1}}{\mathcal F}(g), \ \ {\rm and} \ \ 
vol_{g_{0}}\Omega  \leq  1. 
\end{equation}

 Any compact set $K \subset \Omega_{0}$ embeds in $M$, and if $K$ is 
sufficiently large, the complement $M \setminus K$ carries an 
F-structure, metrically on the complement $B_{i}$ of finitely many balls of 
arbitrarily small radius, as in {\rm (II)}. 

 In both cases {\rm (I)} and {\rm (III)}, the metric $g_{0}$ satisfies the Euler-Lagrange 
equation
\begin{equation} \label{e1.6)}
\nabla{\mathcal F}  = 0. 
\end{equation}
\end{theorem}

  A similar but slightly weaker result holds for the conformally invariant 
functionals ${\mathcal W}^{2}$ or ${\mathcal W}_{\pm}^{2}$, cf. Theorem 4.11. 

 These results give a general framework in which to study the existence 
of minimizers of one of the curvature functionals in (1.1) and, in 
situations where such metrics don't exist on the manifold $M$, a 
framework to try to understand what the obstructions to existence might 
be. Note that the general structure given by Theorem 1.1, and its analogue for 
the conformally invariant functionals, is the same for all functionals 
${\mathcal F}$ in (1.1). 

  The results above also apply to the moduli spaces of minimizers, (or critical 
points) of ${\mathcal F}$, and in this context generalize recent results in [8], [45], 
cf. Theorem 4.15. A number of questions related to Theorem 1.1 are raised in \S 4, 
the most important being to what extent the domain $\Omega$ is topologically 
essential in $M$, analogous to the Thurston decomposition in dimension 3. 

\medskip

 We point out one particular consequence of the proof of Theorem 1.1 here, 
related to an open question of Gromov and work of Rong [37]; again see 
Theorem 4.18 for more details.

\begin{theorem} \label{t 1.2. }
  There is an $\varepsilon_{0} > 0$, such that if $M$ is a 4-manifold 
admitting a metric with 
\begin{equation} \label{e1.7}
\int_{M}|R|^{2} \leq  \varepsilon_{0}, 
\end{equation}
then $M$ has an F-structure. 
\end{theorem}

  We do not attempt here to give a broad overview of results on 
canonical metrics on 4-manifolds, which would require a much longer article; 
thus many important topics are not discussed at all. Some important omissions 
include the existence of canonical K\"ahler-Einstein metrics, where a great 
deal more is known based on Yau's solution of the Calabi conjecture [47]. Similarly, 
extremal K\"ahler metrics and twistor theoretic techniques are not addressed. 
In fact, the relations between the canonical metric problem with complex and algebraic 
geometry are not considered, and it would be interesting to see if the conclusions of 
Theorem 1.1 can be strengthened in the context of K\"ahler metrics for instance.

 Finally, all manifolds below are compact, connected and oriented, and 
of dimension 3 or 4, unless otherwise stated. 

  I would like to thank the referee for several useful comments on the manuscript.

\section{3-Manifolds.}
\setcounter{equation}{0}

 In dimension 3, it is natural to define canonical metrics to be the 
metrics of constant curvature, or equivalently, Einstein metrics. Most 
3-manifolds $M$ do not admit an Einstein metric; in fact it is quite easy 
to see that essential spheres and tori obstruct the existence of an 
Einstein metric, (except tori in the very special case of flat 3-manifolds). 
Here, an embedded sphere $S^{2}$ in $M$ is essential if it does not bound a 
3-ball in $M$, while an embedded torus is essential if the embedding induces 
an injection of fundamental groups. So for example, a non-trivial 
connected sum $M_{1}\#M_{2}$, or any circle bundle over a surface with 
infinite fundamental group, (which is not a flat 3-manifold), does not 
carry an Einstein metric. 

 A special case of the Thurston Geometrization Conjecture, (the most 
important case given Thurston's results on the conjecture [43]), is that 
the simplest essential surfaces embedded in $M$, namely spheres and tori, 
are the only obstructions to the existence of an Einstein metric. In fact, 
the conjecture states that a general 3-manifold may be naturally split 
along a suitable collection of such spheres (sphere decomposition) and 
tori (torus decomposition) into pieces, each of which admits a canonical 
geometric structure. A geometric structure is a mild generalization of an 
Einstein metric, namely a complete, locally homogeneous metric. There are eight 
types of geometries; the three of constant curvature and five which are products 
or twisted products of lower dimensional manifolds, (where the uniformization 
theorem for surfaces comes into play).  

\medskip

 To describe the splitting in a bit more detail, the sphere 
decomposition is a decomposition into a connected sum of irreducible 
3-manifolds, and has the form
\begin{equation} \label{e2.1}
M = (K_{1}\#...\# K_{p})\#(L_{1}\#...\# L_{q})\#(\#_{1}^{r}S^{2}\times S^{1}), 
\end{equation}
where the $K$ and $L$ factors are irreducible and of infinite $\pi_{1}$ 
and finite $\pi_{1}$ respectively; irreducible means that every 
embedded $S^{2}$ bounds a 3-ball $B^{3}$ in the manifold. The torus 
decomposition is a splitting of a $K$-factor into a finite collection of 
disjoint open manifolds $K \setminus \overline{\mathcal T}$, where $\overline{\mathcal T}$ 
is a finite collection of disjoint, non-isotopic, essential tori in $K$ such that 
each component of $K \setminus \overline{\mathcal T}$ has no essential tori not 
homotopic to boundary torus in ${\mathcal T}$. 

 Thurston's conjecture is the assertion that each $L$ factor in (2.1) 
is a spherical space form, while each $K$ factor has the form 
\begin{equation} \label{e2.2}
K = H\cup_{{\mathcal T}}G,  
\end{equation}
where $H$ is a finite union of complete hyperbolic 3-manifolds of 
finite volume, and $G$ is a finite union of graph manifolds, with 
${\mathcal T} \subset \overline{\mathcal T}$. Each component of $G$ may be further 
decomposed as a union of circle bundles over surfaces with boundary, (Seifert 
fibered spaces); the resulting toral boundary components then essentially comprise 
$\overline{\mathcal T} \setminus {\mathcal T}$. Thus $G$ is a union of Seifert fibered 
spaces with boundary, glued together by toral automorphisms. The Seifert fibered 
pieces of $G$ carry product or twisted product geometries. There is one exception 
to the rule (2.2), namely when $K$ is a 3-dimensional Sol-manifold, i.e. a finite 
cover of $K$ is a non-trivial torus bundle over a circle. 

 In a remarkable series of papers [33]-[35], Perelman has apparently 
proved the Geometrization Conjecture. His work has gradually gained 
increasing acceptance among experts and it seems likely that full 
acceptance will occur in the near future. In addition to Perelman's 
papers, there are now a number of expositions of his work at various 
levels, and so the proof will not be discussed here; see also the 
general source [25]. 

\medskip

 For later purposes, there is one point worth explaining however. While 
the method, the Ricci flow with surgery, leads to the geometrization of the 
constant curvature (Einstein) factors in (2.2), it does not lead to the geometric 
structures on the graph manifold part $G$, (or the Sol geometry). Instead, the 
geometry of $G$ that emerges is that of collapse along the circle fibers in the 
Seifert fibered spaces and collapse of the toral regions glueing them together, 
(or collapse of toral fibers in Sol manifolds). Thus, the basic configuration in 
the limit is a collection of Einstein metrics, together with a well-defined 
degeneration by collapse of the remaining parts of $M$. 

 The point worth emphasizing here is that although most 3-manifolds do 
not carry Einstein metrics, given Perelman's work one has a precise 
understanding of which 3-manifolds do, and how a general 3-manifold is 
obtained by assembling pieces having such canonical geometries. 

  Finally, the moduli space of Einstein metrics on 3-manifolds is completely 
understood; the spherical space-forms are rigid (Calabi), as are the hyperbolic 
manifolds of finite volume (Mostow and Mostow-Prasad). The moduli spaces of the 
remaining six geometries are basically determined by the moduli of constant 
curvature metrics on the underlying surfaces. 

\medskip

{\bf Application to the Sigma Constant.}

 Since it does not appear to be widely known at this time, we give an 
application of Perelman's work to the Sigma constant, also called the 
Yamabe invariant, of 3-manifolds. Thus, let ${\mathcal S}$ denote the 
Einstein-Hilbert action restricted to the space ${\mathbb M}_{1}$ of unit 
volume metrics on a given 3-manifold $M$; 
\begin{equation} \label{e2.3}
{\mathcal S} (g) = \int_{M}s_{g}dV_{g}, 
\end{equation}
where $s_{g}$ is the scalar curvature of $g$. ${\mathcal S}$ is bounded 
below in any given conformal class $[g]$ and the invariant $\sigma (M)$ 
is given by 
\begin{equation} \label{e2.4}
\sigma (M) = \sup_{[g]\in{\mathcal C}}\{\inf_{[g]}{\mathcal S}(g)\} = 
\sup_{\gamma\in{\mathcal Y}}s_{\gamma}, 
\end{equation}
where ${\mathcal C}$ is the space of conformal classes and ${\mathcal Y}$ is 
the space of unit volume Yamabe metrics. Now suppose 
\begin{equation} \label{e2.5}
\sigma (M) \leq  0. 
\end{equation}
It follows from classical work of Schoen-Yau or Gromov-Lawson that 
(2.5) occurs if the decomposition (2.1) contains at least one $K$ 
factor; (Perelman's work implies that (2.5) occurs precisely when (2.1) 
contains at least one $K$ factor). 

 We show that Perelman's work implies that when $\sigma (M) \leq 0$, 
$\sigma (M)$ is determined by the volume of the hyperbolic part of $M$, 
in that 
\begin{equation} \label{e2.6}
|\sigma (M)| = 6(vol_{-1}H)^{2/3}, 
\end{equation}
where $vol_{-1}H$ is the volume of $H$ with respect to the metric of 
constant curvature -1. In particular, the graph manifold part $G$ and 
the positive parts $S^{3}/\Gamma$, $S^{2}\times S^{1}$ if any, are invisible 
to $\sigma (M)$. Perelman's work answers affirmatively a conjecture of 
Schoen in [39], and its generalization in [6]. 

 To prove (2.6), consider the quantity 
\begin{equation} \label{e2.7}
S_{-}(M) = \sup_{g\in{\mathbb M}}\{s_{min}v^{2/3}(g)\}, 
\end{equation}
where the sup is taken over the space ${\mathbb M}$ of all metrics on 
$M$, $s_{min}(g) = \min_{M}s_{g}$ and $v$ is the volume of $(M, g)$. The 
product in (2.7) is scale invariant. It is easy to see that when 
$\sigma (M) \leq 0$, then 
\begin{equation} \label{e2.8}
S_{-}(M) = \sigma (M). 
\end{equation}
Namely, since Yamabe metrics are of constant scalar curvature, one has 
$S_{-}(M) \geq  \sigma (M)$. On the other hand, given any $g$, let 
$\widetilde g = u^{4}g$ be a Yamabe metric of the same volume as $g$ in 
$[g]$, so that $u$ satisfies the Yamabe equation
\begin{equation} \label{e2.9}
u^{5}\widetilde s = -8\Delta u + su. 
\end{equation}
When $\widetilde s \leq 0$, the maximum principle implies that 
$\widetilde s \geq  s_{min}$, (since $\max u \geq 1$). This proves 
(2.8), and so (2.6) follows from
\begin{equation} \label{e2.10}
|S_{-}(M)| = 6(vol_{-1}H)^{2/3}. 
\end{equation}

 To prove (2.10), suppose first that $M$ is irreducible, so the sphere 
decomposition (2.1) is trivial, ($M = K$). Then 
\begin{equation} \label{e2.11}
M = H \cup_{{\mathcal T}}G, 
\end{equation}
where the union is along incompressible tori ${\mathcal T}$. Now it is 
easy to construct a metric $g_{\varepsilon}$ on $M$ such that 
\begin{equation} \label{e2.12}
s_{min}v^{2/3}(g_{\varepsilon}) \geq  -6(vol_{-1}H)^{2/3} - \varepsilon , 
\end{equation}
for any given $\varepsilon > 0$. This can be done ``by hand'', by taking a 
truncation of the hyperbolic metric on $H$, joined with a highly 
collapsed metric on $G$; one can easily construct such metrics on $G$ 
with $s \geq -6$, $vol \leq \varepsilon$, for any given $\varepsilon > 0$, 
and which smoothly glue onto the hyperbolic cusps sufficiently far down 
the cusps, cf. [5], [7] for further details. (If $H = \emptyset$, then this 
already implies (2.10)). Thus one has  
\begin{equation} \label{e2.13}
S_{-}(M) \geq  -6(vol_{-1}H)^{2/3}= -\tfrac{3}{2}(vol_{-1/4}H)^{2/3}. 
\end{equation}
Now suppose there is a metric $g_{0}$ on $M$ such that $S_{-}(g_{0}) > 
-\frac{3}{2}(vol_{-1/4}H)^{2/3}$. Then start the Ricci flow on $M$ with 
initial metric $g_{0}$. Perelman's work implies that the Ricci flow 
with surgery $g_{t}$ exists for all time, and that the scale invariant 
quantity $S_{-}(g_{t})$ is monotone non-decreasing in $t$, since 
$S_{-}(g_{t}) \leq 0$ for all $t$, cf. [34]. Hence as $t \rightarrow  
\infty$, $S_{-}(g_{t}) \rightarrow  \widetilde S > 
-\frac{3}{2}(vol_{-1/4}H)^{2/3}$, with $\widetilde S \leq 0$. 

 On the other hand, as Perelman shows, the rescaled metrics 
$\widetilde g_{t} = t^{-1}g_{t}$ have the property that 
$s_{min}(\widetilde g_{t}) \rightarrow  -\frac{3}{2}$ as $t \rightarrow  
\infty$, cf. [34]. Now the decomposition (2.11) is unique up to 
isotopy, (cf. [7]), and the metrics $\widetilde g_{t}$, when restricted to 
compact subsets of $H$, converge to the hyperbolic metric with curvature $-1/4$. 
Thus, one must have $\widetilde V = 
\liminf_{t\rightarrow\infty}vol(\widetilde g_{t}) \geq  vol_{-1/4}H$. 
Hence, $\widetilde S = 
\limsup_{t\rightarrow\infty}s_{min}(\widetilde g_{t})vol(\widetilde g_{t})^{2/3} 
\leq  -\frac{3}{2}(vol_{-1/4}H)^{2/3}$, which gives a contradiction. 

 If $M$ is not irreducible, then $M$ is a connected sum of positive 
factors $S^{3}/\Gamma$, $S^{2}\times S^{1}$ and non-positive irreducible 
factors $K_{i}$. The work above shows that (2.10) holds on each 
$K_{i}$. One can perform the connected sum surgery by hand to increase 
$s_{min}$ pointwise and with arbitrarily small change to the volume, 
cf. again [5], [7], so that (2.12) holds for general $M$. One may then 
apply exactly the same argument as before to prove that (2.10) holds, when 
$\sigma (M) \leq 0$. (The Ricci flow with surgery performs the sphere 
decomposition (2.1), and in particular disconnects the factors in (2.1) in 
finite time, while the $K$ factors persist for infinite time). 
{\endproof}

\medskip

 In contrast, no applications of Perelman's ideas have yet been found 
to determine the Sigma constant of the positive 3-manifolds, i.e. 
$S^{3}/\Gamma$; cf. [1], [14] for some recent progress on this 
problem. 

\medskip

 Observe that (2.10) shows that if $(M, g)$ is any closed Riemannian 
3-manifold with $\sigma (M) \leq 0$, then
\begin{equation} \label{e2.14}
s_{g} \geq  -6 \Rightarrow  vol_{g}M \geq  vol_{-1}H, 
\end{equation}
where $H$ is the hyperbolic part of $M$. This gives a very strong 
generalization of results of [13] in dimension 3, and extends their 
results from Ricci curvature to scalar curvature. 

 In fact, (2.14) can easily be generalized somewhat further. Let
\begin{equation} \label{e2.15}
{\mathcal S}_{-}^{3/2}(g) = \int|\min(s_{g}, 0)|^{3/2}dV_{g}. 
\end{equation}
Then it is easy to see that 
\begin{equation} \label{e2.16}
|\sigma (M)|^{3/2} = \inf_{g\in{\mathbb M}_{1}}{\mathcal S}_{-}^{3/2}(g). 
\end{equation}
Namely, the definition (2.4) gives immediately $|\sigma (M)|^{3/2} \geq  
\inf_{g\in{\mathbb M}_{1}}{\mathcal S}_{-}^{3/2}(g)$. On the other hand, given 
any $g \in  {\mathbb M}_{1}$, let $\gamma$ be a unit volume Yamabe metric 
in $[g]$. Setting $g = u^{4}\gamma$, as in (2.9) one has $u^{5}s_{g} = 
-8\Delta u + s_{\gamma}u$. Since $s_{\gamma}$ is a non-positive constant, 
simple calculations give 
$$|s_{\gamma}| = -\int s_{\gamma}dV_{\gamma} = -\int 
s_{g}u^{4}dV_{\gamma} - 8\int u^{-1}\Delta udV_{\gamma} = -\int 
s_{g}u^{4}dV_{\gamma} - 8\int|d\log u|^{2}dV_{\gamma}  $$
$$\leq  \int|\min(s_{g}, 0)|u^{4}dV_{\gamma} \leq  
[\int|\min(s_{g},0)|^{3/2}u^{6}dV_{\gamma}]^{2/3} = ({\mathcal S}^{3/2}(g))^{2/3}.  
$$
This gives $\inf {\mathcal S}^{3/2}(g) \geq |\sigma (M)|^{3/2}$, and so 
(2.16). Hence, (2.14) generalizes to
\begin{equation} \label{e2.17}
{\mathcal S}_{-}^{3/2}(g) \geq  6^{3/2}vol_{-1}H. 
\end{equation}

 In fact, (2.17) reflects the behavior of metrics minimizing the 
functional ${\mathcal S}_{-}^{3/2}$, (or stronger functionals such as 
${\mathcal S}^{2}$) on a given 3-manifold with $\sigma (M) \leq 0$. Thus, 
one may find minimizing sequences $\{g_{i}\}$ for ${\mathcal S}_{-}^{3/2}$ 
which crush essential 2-spheres in $M$ to points, according to the sphere 
decomposition, and on each non-positive $K$-factor, converge to the complete 
hyperbolic metric on $H$, while collapsing the graph manifold part $G$ with 
uniformly bounded curvature. The positive parts $S^{3}/\Gamma$ and 
$S^{2}\times S^{1}$ are invisible to ${\mathcal S}_{-}^{3/2}$. Thus, one can 
construct minimizing sequences for ${\mathcal S}_{-}^{3/2}$ which give a geometric 
decomposition of $M$, equivalent to the Thurston decomposition, (cf. [7] 
for further details). 

 This will be the main point of view in the analysis to follow in 4-dimensions.

\section{4-Manifolds: Moduli Spaces.}
\setcounter{equation}{0}

 On 4-manifolds, it is less clear what a canonical metric should be. As 
discussed in \S 1, we will take the point of view of variational 
problems on the space of metrics ${\mathbb M}$ on a given 4-manifold $M$ and 
define such a metric to be a minimizer, (or possibly a critical point), 
of one of the curvature functionals ${\mathcal F}$ in (1.1), i.e. 
\begin{equation} \label{e3.1}
{\mathcal R}^{2}, \ {\mathcal W}^{2}, \ {\mathcal W}_{+}^{2}, \ {\mathcal W}_{-}^{2}, \ {\mathcal R}ic^{2}, 
\end{equation}
or the much weaker scalar curvature analogues,
\begin{equation} \label{e3.2}
{\mathcal S}^{2}, \ -{\mathcal S}|_{{\mathcal Y}}. 
\end{equation}

 These functionals are all bounded below, and so in principle one can 
use direct methods in the calculus of variations to study the existence 
and properties of minimizers. 

\medskip

 The basic problem is to understand the existence and moduli spaces 
of such metrics on a given manifold $M$. However, just as in dimension 3 
as discussed in \S 2, one cannot expect an arbitrary 4-manifold admits 
a smooth metric minimizing one of the functionals ${\mathcal F}$. In fact, 
the situation in dimension 4 is much more complicated than 
that in 3 dimensions. While numerous obstructions to the existence 
of minimizers of a given ${\mathcal F}$ are known, cf. [28], [30] and further 
references therein, there is no general conjecture as to what an exact and complete 
set of obstructions is, i.e. there is currently no analog of the 
Thurston geometrization conjecture. 

 Nevertheless, it is natural to try to find a geometric decomposition 
of $M$ with respect to one of these functionals. Thus, as discussed at the 
end of \S 2, one can try to see if minimizing sequences decompose the manifold 
into pieces, (analogous to (2.1) or (2.2)), on some of which they converge to 
smooth limits and others on which they degenerate in a well-defined way. 

 The single most important fact allowing one to develop such a theory 
on the existence, or the structure of moduli spaces of such functionals, 
is Chern's generalization of the Gauss-Bonnet theorem; in dimension 4, 
this is
\begin{equation} \label{e3.3}
\frac{1}{8\pi^{2}}\int\{|R|^{2} - |z|^{2}\}dV = \chi (M), 
\end{equation}
where $z = Ric - \frac{s}{4}g$ is the tracefree Ricci curvature. The 
expression (3.3) is equivalent to
\begin{equation} \label{e3.4}
\frac{1}{8\pi^{2}}\int\{|W|^{2} - \frac{1}{2}|z|^{2} + \frac{1}{24}s^{2}\}dV 
= \chi (M). 
\end{equation}
This gives one $L^{2}$ control of the full curvature $R$ (or $W$) in 
terms of $L^{2}$ control of $Ric$. Chern-Weil theory and the signature 
theorem also give the relation
\begin{equation} \label{e3.5}
\frac{1}{12\pi^{2}}\int\{|W_{+}|^{2} - |W_{-}|^{2}\}dV = \tau (M), 
\end{equation}
where $\tau (M)$ is the signature of $M$. Combining (3.4) and (3.5) gives
\begin{equation} \label{e3.6}
\frac{1}{2\pi^{2}}\int\{|W_{+}|^{2} - \frac{1}{4}|z|^{2} + 
\frac{1}{48}s^{2}\}dV = 2\chi (M) + 3\tau (M). 
\end{equation}

 The functionals (3.1) are all scale-invariant in dimension 4. In the 
following, we will always work on the space ${\mathbb M}_{1}$ of unit 
volume metrics on $M$, unless stated otherwise.  

\medskip

 In this section, we study the structure of the moduli spaces of minimizers 
or critical points of the functionals in (3.1). This serves as an introduction 
as to what one can expect for existence results, which are discussed in \S 4. 

\medskip

 (A). Einstein Moduli Spaces.

\medskip

 We begin with the case of Einstein metrics, which are critical points 
of all the functionals in (3.1) and (3.2). Let ${\mathcal M}  = {\mathcal M}_{E}$ 
denote the moduli space of unit volume Einstein metrics on $M$. 
Since for instance the functional ${\mathcal S}^{2}$ is critical on ${\mathcal M}$ 
and Einstein metrics have constant scalar curvature, the scalar 
curvature $s_{g}: {\mathcal M}  \rightarrow  {\mathbb R}$ is constant on 
components of ${\mathcal M}$. By (3.3), ${\mathcal R}^{2}$ is constant on all of 
${\mathcal M}$, while by (3.4), ${\mathcal W}^{2}$ is again constant on components 
of ${\mathcal M}$. 

\medskip

 The first general result on the structure of the moduli space ${\mathcal M}$ of 
unit volume Einstein metrics on a given 4-manifold $M$ was obtained in [2], [10], 
[32]; a partial result along these lines was also obtained in [44] in the special 
case of K\"ahler-Einstein metrics with $c_{1} > 0$. Overall, the picture resembles 
somewhat Uhlenbeck's results on the moduli space of self-dual Yang-Mills fields. 

 To describe this, an Einstein orbifold $(V, g)$ associated to $M$ 
is defined to be a (4-dimensional) orbifold, with a finite number of singular 
points $\{q_{i}\}$, each having a neighborhood homeomorphic to the cone 
$C(S^{3}/\Gamma)$, where $\Gamma \neq \{e\}$ is a finite subgroup of $SO(4)$. 
Let $V_{0} = V\setminus \cup q_{k}$ be the regular set of $V$. Then $g$ 
is a smooth Einstein metric on $V_{0}$, which extends smoothly over 
$\{q_{k}\}$ in local finite covers. The manifold $M$ is a resolution 
of $V$ in the sense that there is a continuous surjection 
$\pi : M \rightarrow  V$ such that $\pi|_{\pi^{-1}(V_{0})}: 
\pi^{-1}(V_{0}) \rightarrow  V_{0}$ is a diffeomorphism onto $V_{0}$. 
In particular, $V$ is compact.  

 Then the result is that the completion $\hat {\mathcal M}$ of ${\mathcal M}$ in the 
Gromov-Hausdorff topology consists of ${\mathcal M}$ together with unit 
volume Einstein orbifold metrics associated to $M$. Moreover, the 
completion is locally compact, in that any sequence $g_{i}$ of unit 
volume Einstein metrics on $M$, bounded in the Gromov-Hausdorff topology, 
has a subsequence converging to an Einstein orbifold associated to $M$. 

 In analogy to the Uhlenbeck completion of the moduli space of 
Yang-Mills instantons, orbifold singularities arise from the bubbling 
off of gravitational instantons, i.e. complete non-flat Ricci-flat metrics 
$(N, g_{\infty})$ which are ALE, (asymptoticaly locally Euclidean), in that 
the metric $g_{\infty}$ is asymptotic to a flat cone $C(S^{3}/\Gamma)$ 
at infinity. There is at least one such ALE space associated to each 
singularity; however, in general there may be finitely many such spaces, 
arising at different blow-up scales. All such blow-up limits $N$ are 
topologically embedded in $M$ and a simple Mayer-Vietoris argument shows 
that most of the rational homology of any such $N$ injects in the homology 
of $M$, in that 
\begin{equation} \label{e3.7}
0 \rightarrow H_{k}(N, {\mathbb R}) \rightarrow H_{k}(M, {\mathbb R}),
\end{equation}
for $k = 1, 2$. 

  Such ALE spaces $(N, g_{\infty})$ have nontrival $2^{\rm nd}$ 
Betti number, and $V$ is obtained from $M$ by collapsing essential 
cycles in $H_{2}(M, {\mathbb R})$ to points. In particular, if $b_{2}(M) = 0$, 
then there are no orbifold singularities and so $V = M$. This is the 
case for example $M$ is a rational homology sphere, (with any differentiable 
structure). Also, the proof of the smoothness of $g$ across the orbifold 
singularities in a local uniformization follows the lines of proof of Uhlenbeck's 
removable singularity theorem, cf. [10]. Finally, since each ALE space 
$(N, g_{\infty})$ has a definite amount of curvature in $L^{2}$, one has 
a uniform bound on the number of orbifold singularities depending only on 
$\chi (M)$, by the Chern-Gauss-Bonnet theorem (3.3).

 The completion in the Gromov-Hausdorff topology is equivalent to the 
completion with respect to a diameter bound, so that all metrics 
$g \in \hat{\mathcal M}$ satisfy
\begin{equation} \label{e3.8}
vol_{g}M = 1, \ \ diam_{g}M \leq  D,  
\end{equation}
for some $D = D(g) < \infty$. 

  However, there is a marked difference compared with the Uhlenbeck completion. 
Namely, the frontier $\hat {\mathcal M} \setminus {\mathcal M}$ should perhaps not be 
considered as a boundary, but instead as a filling in of missing pieces in ${\mathcal M}$. 
For example, in the case of K3 surfaces, the moduli space has dimension 57, and 
the frontier consists of subvarieties of codimension 3 in $\hat{\mathcal M}$, cf. [12]; 
it does not form a natural boundary as a ``wall'' past which $\hat {\mathcal M}$ cannot 
be continued. This is the case in all known examples, although it is unknown if 
this holds in general. 

 A main point of the Uhlenbeck completion is that the completion is 
compact. Consider the components of ${\mathcal M}$ for which 
\begin{equation} \label{e3.9}
s_{g} \geq  s_{0} >  0. 
\end{equation}
Myers' theorem then implies (3.8) holds, (with $D = D(s_{0}))$, and so 
the completion $\hat{\mathcal M}$ of this part of ${\mathcal M}$ in the Gromov-Hausdorff 
topology is compact. However, this is certainly not the case when $s_{g} \leq 0$, 
(consider for example flat tori or products of hyperbolic surfaces).  Thus, 
one needs to consider what happens when the Gromov-Hausdorff distance 
goes to infinity. 

\medskip

 In [4] a more complete theory of the global behavior of ${\mathcal M}$ 
was developed, which from a broad perspective has a strong resemblance 
with the moduli space of constant curvature metrics on surfaces. (The 
compactness of the part of $\hat {\mathcal M}$ for which (3.9) holds 
corresponds to the compactness of the moduli space of Einstein metrics 
on $S^{2})$. In studying the boundary of Teichm\"uller space or the Riemann 
moduli space, one of the most natural metrics is the Weil-Petersson metric. 
This has a natural generalization to all dimensions, since it is just the 
restriction of the natural $L^{2}$ metric on ${\mathbb M}$ to the moduli space. 
Thus, we consider the completion of ${\mathcal M}$ with respect to the $L^{2}$ 
metric. 

 To describe the results, one needs the following definition. A domain $\Omega$ 
(i.e. an open 4-manifold) weakly embeds in $M$, $\Omega \subset \subset M$, 
if for any compact subdomain $K \subset  \Omega$, there is a smooth embedding 
$F = F_{K}: K \rightarrow M$. The same definition applies if $\Omega$ is an 
orbifold, with the obvious modification that the corresponding part of $M$ 
is a resolution of $K$. 

 The completion $\overline{\mathcal M}$ of ${\mathcal M}$ with respect to the 
$L^{2}$ metric on ${\mathbb M}$ is a complete Hausdorff metric space, 
whose frontier $\partial{\mathcal M}$ consists of two parts: the orbifold 
part $\partial_{o}{\mathcal M}$ and the cusp part $\partial_{c}{\mathcal M}$. 

 (I). $\partial_{o}{\mathcal M}$ consists of compact Einstein orbifolds of 
unit volume associated to $M$. The partial completion 
${\mathcal M}\cup\partial_{o}{\mathcal M}$ is locally compact. This is the 
same as the situation described before. 

 (II). An element in the cusp boundary $\partial_{c}{\mathcal M}$ is given 
by a pair $(\Omega, g)$, where $\Omega$ is a non-empty maximal 
orbifold domain $\Omega$ weakly embedded in $M$. The domain $\Omega$ consists 
of a finite number of components $\Omega_{k}$ called cusps, each with a bounded 
number, (possibly zero), of orbifold singularities. The metric $g$ is a complete 
Einstein metric on $\Omega$, with 
\begin{equation} \label{e3.10}
vol_{g}\Omega  = 1 ,
\end{equation}
and outside a compact set $K$, $\Omega$ carries an $F$-structure along 
which $g$ collapses with locally bounded curvature as one goes to 
infinity in $\Omega$; thus as $x \rightarrow  \infty$ in $\Omega$, 
\begin{equation} \label{e3.11}
inj(x) \rightarrow 0 \ \ {\rm and} \ \ (|R|inj)^{2}(x) \rightarrow 0, 
\end{equation}
where $inj(x)$ is the injectivity radius at $x$. 

 To describe the behavior of the region $M\setminus K$, let $g_{i}$ 
be a sequence in ${\mathcal M}$ with $g_{i} \rightarrow g$ in the $L^{2}$ 
metric. Then $M \setminus K$ also carries an $F$-structure on the 
complement of a finite number of arbitrarily small balls. Thus, there 
exists a finite collection of points $z_{j}\in M$, and a sequence 
$\varepsilon_{i} \rightarrow 0$ such that outside 
$B_{z_{j}}(\varepsilon_{i})$, $M \setminus K$ has an $F$-structure. If 
one chooses an exhaustion $K_{i}$ of $\Omega$, then 
$M\setminus K_{i}$ collapses everywhere, and collapses with locally 
bounded curvature away from the singular points $\{z_{j}\}$. 

 The convergence in (I) is also in the Gromov-Hausdorff topology, while 
that in (II) is also in the pointed Gromov-Hausdorff topology, for a suitable 
collection of base points. Further, cusps, i.e. Case (II), can occur only on 
the components of ${\mathcal M}$ for which there is a constant $s_{0}$ such that
\begin{equation} \label{e3.12}
s_{g} \leq  s_{0} < 0. 
\end{equation}

 The $L^{2}$ completion $\overline{\mathcal M}$ is not compact in general, 
as seen explicitly in the cases of flat metrics on tori, or Ricci-flat 
metrics on $K3$. On the $K3$ surface, the $L^{2}$ metric on $\overline{\mathcal M}$ 
is the complete metric of finite volume on the non-compact locally symmetric 
space $\Gamma \setminus SO(3,19)/(SO(3)\times SO(19))$, cf. [4], [12], so 
that $\overline{\mathcal M} = {\mathcal M} \cup \partial_{o}{\mathcal M}$. 
(Of course if (3.9) holds, then $\overline{\mathcal M}$ is 
compact). The behavior of $\overline{\mathcal M}$ at infinity is described 
as follows:

 (III). Suppose $g_{i}$ is a divergent sequence in $\overline{\mathcal M}$ 
such that 
\begin{equation} \label{e3.13}
s_{g_{i}} \rightarrow  0, \ \ {\rm as} \ \ i \rightarrow  \infty . 
\end{equation}
Then $\{g_{i}\}$ collapses everywhere with locally bounded curvature, 
i.e. (3.11) holds, metrically on the complement of finitely many singular 
points $\{z_{j}\}$. (Thus, $\Omega = \emptyset$ in the context of (II)). 

 Suppose instead that $g_{i}$ is a divergent sequence in $\overline{\mathcal M}$ 
such that 
\begin{equation} \label{e3.14}
s_{g_{i}} \leq  s_{0} <  0, \ \ {\rm as} \ \ i \rightarrow  \infty . 
\end{equation}
Then $\{g_{i}\}$ either has the same behavior as in (III) or (II), where 
$\Omega$ may instead have possibly infinitely many components, (of 
total volume at most 1). 

 Recently, building on the work in [4], Cheeger-Tian [16] have 
improved the statement above, and proved:

 (IV). Suppose instead that $g_{i}$ is a divergent sequence in 
$\overline{\mathcal M}$ such that 
$$s_{g_{i}} \leq  s_{0} <  0, \ \ {\rm as} \ \ i \rightarrow  \infty . $$
Then $\{g_{i}\}$ has the same behavior as in (II). Moreover, the collapse with 
locally bounded curvature is actually collapse with uniformly bounded 
curvature: $|R| \leq \Lambda$, for some $\Lambda = \Lambda(M) < \infty$, 
away from the singular points.  

 To complete the analogy with the case of surfaces, it is natural to 
conjecture, (cf. [4]), that in fact Case (IV) does not occur, i.e. when 
(3.14) holds, $\overline{\mathcal M}$ is compact. We recall here that the 
completion of the moduli space of hyperbolic metrics on a surface with 
respect to the Weil-Petersson $(L^{2})$ metric is compact, and agrees 
with the Deligne-Mumford compactification. 

  There are many open questions regarding the structure of ${\mathcal M}_{E}$ 
that remain unanswered. Among the most basic are the following: 

(i). Does ${\mathcal M}$ or $\overline{\mathcal M}$ have finitely many components? 
This is open even for the portion of ${\mathcal M}$ satisfying (3.9) where 
$\overline{\mathcal M}$ is compact. Apriori, one could have a sequence of 
metrics $g_{i}$ in distinct components of ${\mathcal M}$ which converge to 
an orbifold metric $(V, g) \in \overline{\mathcal M}$ in a limit component. 

(ii). If $g_{i}$ is a sequence in ${\mathcal M}$ converging to $(V, g) \in 
\overline{\mathcal M}$ as above, does there exist a curve $\gamma(t)$ in 
$\overline{\mathcal M}$, with $\gamma(0) = g$ and $\gamma(t_{i}) = g_{i}$, 
for some sequence $t_{i} \rightarrow 0$?

(iii). The occurence of orbifold singularities is closely related to the topology 
of $M$, via (3.7). If it is possible topologically for orbifold singularities 
to occur, do they in fact occur in $\overline{\mathcal M}$?

(iv). Is there any relation between the topology of $\Omega$ and the topology of 
$M$, i.e. is $\Omega$ in any way topologically essential in $M$?

(v). What can be said about the structure of the singularities in the collapsing 
situation of Case (III), or in the complement of $\Omega$ in Case (II)? A concrete 
example of collapsing metrics on the $K3$ surface is described in detail in [23]; 
the singularities here are modeled on the Ooguri-Vafa metrics, which are periodic 
versions of the Taub-NUT metrics. 

\medskip

 (B) General Moduli spaces. 

\medskip

 The first result (I) discussed above, on the completion of the moduli 
space of Einstein metrics on $M^{4}$ with respect to the 
Gromov-Hausdorff topology, has recently been generalized to the 
functionals in (3.1). Thus, let ${\mathcal F}$ denote any of the 
functionals in (3.1), and consider the moduli space ${\mathcal M}_{{\mathcal F}}$ 
of critical points of ${\mathcal F}$ of unit volume. As before, ${\mathcal F}$ is 
constant on each component of ${\mathcal M}_{\mathcal F}$. 

 First, the definition of orbifold needs to enlarged somewhat, in that 
one allows a neighborhood of a singular point $q$ to be given by a 
finite collection of cones $C(S^{3}/\Gamma_{j})$, with vertex $q$, (not 
just a single cone). For emphasis, sometimes such orbifolds will be called 
reducible, while the orbifolds as previously defined will be called irreducible. 
In particular, an orbifold $V$ is reducible if and only if the regular set 
$V_{0}$ has more than one component near each singular point. Further, some 
or all of the finite groups $\Gamma_{j}$ may be trivial, corresponding to cones 
on $S^{3}$ and so 4-balls. Also, from now on, the metric $g$ on $V$ is only 
asserted to be $C^{0}$ across each singularity $q$ in local uniformizations 
of the cones. 

 Choose $\nu_{0} > 0$ and $\Lambda < \infty$, and let 
${\mathcal M}_{{\mathcal F}}(\nu_{0}, \Lambda)$ denote the portion of 
${\mathcal M}_{{\mathcal F}}$ consisting of all unit volume metrics $g$ such that, 
for $r \leq 1$,
\begin{equation} \label{e3.15}
vol B_{x}(r) \geq  \nu_{0}r^{4}, \ \ {\rm and} \ \ \int |R|^{2} \leq \Lambda ,
\end{equation}
where $B_{x}(r)$ is the geodesic $r$-ball about $x$ in $(M, g)$. Then for 
any $\nu_{o} > 0$ and $\Lambda < \infty$, the closure of 
${\mathcal M}_{{\mathcal F}}(\nu_{0}, \Lambda)$ in the Gromov-Hausdorff topology 
consists of $\hat {\mathcal M}_{{\mathcal F}} = 
{\mathcal M}_{{\mathcal F}}\cup  {\mathcal M}_{o}$, where ${\mathcal M}_{o}$ consists 
of orbifold singular metrics $(V, g)$, as in (I), with the modifications in 
the definition of orbifolds discussed above; in particular, the orbifold $V$ 
may be reducible. Here, in the case of the conformally invariant functionals 
${\mathcal W}^{2}$ and ${\mathcal W}_{\pm}^{2}$, the representative metrics $g\in [g]$ 
are assumed to be unit volume Yamabe metrics. This result was proved in [8] 
and [45]. Under certain conditions, e.g. the presence of a uniform Sobolev 
inequality in addition to (3.15), the orbifolds are irreducible. 

\medskip

 Of course the condition (3.15) rules out any collapse behavior. In \S 4, we 
discuss analogs of the general degeneration for Einstein metrics for the moduli 
spaces ${\mathcal M}_{{\mathcal F}}$, (cf. Theorem 4.15). 

  At this point, it is useful to consider the following simple example, which 
illustrates some strong differences between the functionals ${\mathcal R}ic^{2}$ or 
${\mathcal R}^{2}$, and ${\mathcal W}^{2}$. 
\begin{example} \label{ex3.1}
{\rm Let $M = \Sigma \times S^{1}$, where $\Sigma$ is a compact hyperbolic 3-manifold and 
let $g_{\mu}$ be a product metric on $M$ of the form 
\begin{equation} \label{e3.16}
g_{\mu} = \mu^{2}g_{-1} + \mu^{-6}g_{S^{1}(1)},
\end{equation}
where $g_{-1}$ is the hyperbolic metric on $\Sigma$. The metrics $g_{\mu}$ are 
conformally flat, $W = 0$, and so give a curve in the moduli space ${\mathcal M}_{\mathcal W}$ of 
minimizers of ${\mathcal W}^{2}$ on $M$. 

  The scalar curvature of $g_{\mu}$ is given by $s_{g_{\mu}} = -6\mu^{-2}$ while the 
volume is $vol_{g_{\mu}}M = 2\pi vol_{g_{-1}}\Sigma$. As $\mu \rightarrow 0$, one thus 
has 
\begin{equation} \label{e3.17}
{\mathcal S}^{2}(g_{\mu}) = \int_{M}s_{g_{\mu}}^{2}dV_{g_{\mu}} \rightarrow \infty.
\end{equation}
In particular, the $L^{2}$ norm of the curvature $R$ of $g_{\mu}$ diverges to infinity, 
in strong contrast to the situation of metrics with bounds on ${\mathcal R}ic^{2}$ or 
${\mathcal Z}^{2}$. 

  As $\mu \rightarrow 0$, the factor $\Sigma$ collapses to 0 volume, causing the curvature 
to blow up, while the $S^{1}$ factor expands, to preserve a fixed volume. This behavior 
is completely different from the kind of limits one sees in the moduli space of Einstein metrics 
in (I)-(IV) above. In the opposite direction where $\mu \rightarrow \infty$, the hyperbolic 
factor expands, (to a flat metric), while the circle $S^{1}$ collapses; this behavior is of 
course consistent with the collapse behavior discussed above. 

  One may also perform the same construction with $(\Sigma, g_{-1})$ replaced by 
$(S^{3}, g_{+1})$, preserving conformal flatness. The estimate (3.17) holds, 
where now the scalar curvature diverges to $+\infty$ instead of $-\infty$. Again the 
metrics $g_{\mu}$ behave badly as $\mu \rightarrow 0$, and collapse in a quite different 
manner than the Einstein case. 

  Here however, the divergence of (3.17) and $g_{\mu}$ is due to only to a bad choice of 
gauge for the conformal class. While the metrics $g_{\mu}$ are Yamabe metrics in the negative 
case, they are not Yamabe in the positive case, since for instance there is a uniform upper 
bound on the scalar curvature of Yamabe metrics with a fixed volume. In a Yamabe gauge, 
the conformal classes $(S^{3}\times S^{1}, [g_{\mu}])$ have a uniform bound on the $L^{2}$ 
norm of curvature, and converge to the round metric on $S^{4}$, with two antipodal points 
identified, giving rise to an orbifold singularity consisting of two cones on $S^{3}$ joined 
at the vertex, cf. [39], [8]. }
\end{example}

\section{4-Manifolds: Existence Issues.}
\setcounter{equation}{0}

 One would like to extend the results above on the moduli spaces toward 
an existence theory for metrics minimizing one of the functionals 
${\mathcal F} $ in (3.1), (or (3.2)). From the point of view of direct methods 
in the calculus of variations, this requires understanding the limiting 
behavior of minimizing sequences of metrics $g_{i}$ for ${\mathcal F}$ on a given 
4-manifold $M$. 

 Consider for instance ${\mathcal R}^{2}$, with absolute minimum realized by 
Einstein metrics, (if such exist), by the Chern-Gauss-Bonnet theorem (3.3). 
Any sequence $g_{i}\in{\mathcal M}_{E}$ is trivially a minimizing sequence for 
${\mathcal R}^{2}$, and so general minimizing sequences will exhibit at least the 
degenerations described in (I)-(IV) of \S 3. In other words, a general existence 
theory of minimizers of ${\mathcal F}$ must include a full description of the moduli 
space ${\mathcal M}_{{\mathcal F}}$. 

 However, a general minimizing sequence $\{g_{i}\}$ say for ${\mathcal 
R}^{2}$ will not satisfy any particular PDE, (as is the case on the 
moduli space). At best, the bound on ${\mathcal R}^{2}$ gives $L^{2,2}$ 
control of the metrics $g_{i}$ in some local coordinate charts. 
This is troublesome in dimension 4, since an $L^{2,2}$ bound does not 
give a pointwise $(L^{\infty})$ bound. With only such weak control, it 
does not seem possible to say anything reasonable about some limiting 
behavior of $\{g_{i}\}$. 

 A different approach to a rather general existence theory was 
introduced and developed by Taubes, both for the existence of self-dual 
Yang-Mills instantons [40], and self-dual metrics, (which of course are 
minimizers of ${\mathcal W}_{-}^{2}$), [41], cf. also [42]. The idea here is 
to glue together exact solutions on pieces of a manifold to obtain a global 
approximate solution, and show the approximation can be perturbed into a 
global exact solution. Ingeniously, to accomplish this Taubes introduces a 
slightly stronger norm than $L^{2,2}$ norm, which does embed in $L^{\infty}$. 
Further comments on this approach will follow later in \S 4. 

 The problem that $L^{2,2}$ does not embed in $L^{\infty}$ can be overcome if 
instead one minimizes the $L^{2p}$ norm of the curvature tensor, for any $p > 1$. 
In fact, there is a general ``convergence theorem'' in this context, which shows 
that many of the main features of the results in \S 3 for the moduli space 
${\mathcal M}_{E}$ hold for metrics with just $L^{2p}$ bounds on the Ricci 
curvature.

\begin{theorem} \label{t 4.1.}
  (Convergence under $L^{p}$ curvature bounds.) Let $\{g_{i}\}$ be a 
sequence of unit volume metrics on $M$ such that
\begin{equation} \label{e4.1}
\int_{M}|Ric|^{2p} \leq  \Lambda , 
\end{equation}
for some $p > 1$ and $\Lambda < \infty$. 

 Then a subsequence of $\{g_{i}\}$ converges in the pointed 
Gromov-Hausdorff topology to a maximal orbifold domain $\Omega$, 
possibly empty, weakly embedded in $M$, with $L^{2,2p}$ smooth metric $g$ 
on the regular set $\Omega_{0}$. The convergence is in the weak 
$L^{2,2p}$ topology on $\Omega_{0}$ and so in particular $g$ is locally 
$C^{\alpha}$ on $\Omega_{0}$, $\alpha  = 1 - (2-p)/p$. The orbifold 
singularities are all irreducible, with local group $\Gamma \neq \{e\}$, 
and the metric extends $C^{0}$ across orbifold singularities. One has 
\begin{equation} \label{e4.2}
vol_{g}\Omega  \leq  1. 
\end{equation}
Outside a sufficiently large compact set $K \subset \Omega$, $\Omega$ 
carries an $F$-structure along which $g$ collapses with locally bounded 
curvature in $L^{2p}$ on approach to $\partial \Omega$. The complement 
$M \setminus K$ also carries an $F$-structure on the complement of a finite 
number of balls $B_{z_{j}}(\varepsilon_{i})$ with respect to $g_{i}$, where 
$\varepsilon_{i} \rightarrow 0$ as $i \rightarrow \infty$. If $K_{i}$ is an 
exhaustion of $\Omega$, then $M\setminus K_{i}$ collapses 
everywhere, and collapses with locally bounded curvature in $L^{2p}$ away 
from the singular points $\{z_{j}\}$. 
\end{theorem}

 We point out that in general, the metric $g$ will not be complete on 
$\Omega$, and the domain $\Omega$ could have infinitely many components. 
On the other hand, $\Omega$ may be compact, so that $\Omega = M$ or $\Omega$ 
is an orbifold $V$. Further, one may have $\Omega = \emptyset$, in which case 
$\{g_{i}\}$ collapses $M$ along a sequence of F-structures with locally 
bounded curvature in $L^{2p}$ on the complement of a bounded number of 
arbitrarily small balls.

  Also, here and in the following, we will not detail the choice of (the 
collection of) base points used for the pointed Gromov-Hausdorff topology. 
Some choices of base points lead to limits consisting of components of 
$\Omega$, while other choices lead to the collapsed part on the complement 
of $\Omega$. 

\medskip

 For the proof of Theorem 4.1, we first need the following definitions, 
cf. [3], [5]. 

\noindent
{\bf Definitions.}
 Let $(N, g)$ be a Riemannian 4-manifold. Let $B_{x}(r)$ denote the 
geodesic $r$-ball about $x$ in $N$ and for $y\in B_{x}(r)$, let $D_{y}(s) 
= B_{x}(r)\cap B_{y}(s)$, $s \leq r$. The $\mu$-{\it volume radius} of $(N, g)$ 
at $x$ is 
\begin{equation} \label{e4.3}
\nu (x) = \nu_{\mu}(x) = \sup\{r: \frac{vol D_{y}(s)}{s^{4}} \geq  
\mu\}, 
\end{equation}
for all $y\in B_{x}(r)$. The parameter $\mu > 0$ may be freely chosen, 
but from now on we assume $\mu$ is fixed, say $\mu  = 10^{-2}$. 

 The $L^{q}$ {\it curvature radius} $\rho^{(q)}(x)$ is the largest radius 
such that for $y\in B_{x}(\rho^{(q)}(x))$, one has
\begin{equation} \label{e4.4}
\frac{s^{2q-4}}{vol D_{y}(s)}\int_{D_{y}(s)}|R|^{q} \leq  c_{0}. 
\end{equation}
We assume that $q \geq 2$, so that $L^{2,q} \subset C^{\alpha}$, 
$\alpha  = 1 - (4-q)/q$ when $q > 2$. The parameter $c_{0} = c_{0}(q) > 0$ will 
be chosen to be sufficiently small. An important case is $q = 2$, and we set 
$\rho^{(2)} = \rho$. 

 The $L^{2,q}$ harmonic radius $r_{h}^{2,q}$ at $x$ is the radius of 
the largest geodesic ball about $x$ for which there exists a harmonic 
coordinate chart on $B_{x} = B_{x}(r_{h}^{2,q}(x))$ in which the metric 
components $g_{ij}$ satisfy
\begin{equation} \label{e4.5}
e^{-C}\delta_{ij} \leq  g_{ij} \leq  e^{C}\delta_{ij}, \ \ {\rm as \ bilinear 
\ forms}, 
\end{equation}
and
\begin{equation} \label{e4.6}
(r_{h}^{2,q})^{\lambda}||\partial^{2}g_{ij}||_{L^{q}(B_{x})} \leq  C, 
\end{equation}
where $\lambda = (2q-4)/q$; here one assumes $q > 2$. 

\medskip

  The parameter $C$ is a fixed constant, e.g. $C = 1$. The radii $\nu$, 
$\rho^{(q)}$ and $r_{h}^{2,q}$ all scale as distances under rescalings of 
the metric. 

 For $q > 2$, one has $\rho^{(q)}(x) \geq  c \cdot r_{h}^{2,q}(x)$, for a fixed 
numerical constant $c$, (depending only on $c_{0},C$). More importantly, 
there is constant $c_{1} > 0$, depending only on $c_{0}, C$ and a lower bound 
$\nu_{0}$ for $\nu(x)$, such that
\begin{equation} \label{e4.7}
\rho^{(q)}(x) \leq  c_{1}r_{h}^{2,q}(x).  
\end{equation}
Thus, (4.7) holds on scales bounded above by the volume radius. In the following, 
it will be assumed that $c_{0}$ is chosen sufficiently small so that the ball 
$B = B_{x}(\min(\nu (x), \rho^{(q)}(x))$ is diffeomorphic to a ball in a flat 
manifold and the metric is $C$-close to the flat metric on $B$. 

 Finally, in harmonic coordinates, the Ricci curvature is an elliptic 
operator in the metric; in fact
\begin{equation} \label{4.8}
-\tfrac{1}{2}Ric_{ij} = \Delta_{g}g_{ij} + Q_{ij}(g, \partial 
g), \ {\rm where} \  \Delta_{g} = g^{ab}\partial_{a}\partial_{b}.
\end{equation}
Hence, if the metric is controlled in $C^{\alpha}$, for some $\alpha 
> 0$, then an $L^{q}$ bound for $Ric$ implies an $L^{2,q}$ bound for $g$, 
and hence an $L^{q}$ bound for the full curvature $R$. 

\medskip

{\bf Proof of Theorem 4.1.}
 Let $\{g_{i}\}$ be any sequence of unit volume metrics on $M$ 
satisfying (4.1); in the following we will usually write $g$ for any of 
the metrics $g_{i}$. As will be seen, a crucial point is that the bound 
(4.1), together with the Chern-Gauss-Bonnet theorem (3.3) gives a bound
\begin{equation} \label{e4.9}
\int|R|^{2} \leq  \Lambda' , 
\end{equation}
where $\Lambda'$ depends only on $\Lambda$ and $\chi(M)$. 

 Pick any $\nu_{0} > 0$, and let $M^{\nu_{0}}$ be the $\nu_{0}$-thick 
part of $(M, g)$, $(g = g_{i})$, given by
\begin{equation} \label{e4.10}
M^{\nu_{0}} = \{x\in M: \nu (x) >  \nu_{0}\}. 
\end{equation}
Then $M^{\nu_{0}}$ is an open submanifold in $M$, and it follows from 
[3] that is $L^{2,2p}$ orbifold compact, (away from its boundary), 
in that $(M^{\nu_{0}}, g_{i})$ has a subsequence converging to a 
limit orbifold domain $(\Omega_{\nu_{0}}, g)$, embedded in $M$; (more 
precisely, there is a domain in $M$ which is a resolution of 
$(\Omega_{\nu_{0}}, g)$). The metric $g$ is an $L^{2,2p}$ metric, 
$C^{0}$ across the orbifold singularities, and $\Omega_{\nu_{0}}$ has 
finitely many components. There is a uniform bound on the number of 
orbifold singular points on $\Omega_{\nu_{0}}$, depending only on 
$\chi (M)$, for the same reasons as discussed following (3.7); namely 
each orbifold singularity is associated to a finite number of Ricci-flat 
ALE spaces, each of which contributes a definite amount to the $L^{2}$ 
norm of the curvature in (4.9). Similarly, the orbifold $\Omega_{\nu_{0}}$ 
is irreducible, with local group $\Gamma \neq \{e\}$. 

 Let $\nu_{j}$ be a decreasing sequence with $\nu_{j} \rightarrow 0$ 
as $j \rightarrow  \infty$. The argument above applies to each 
$M^{\nu_{j}} \subset M^{\nu_{j+1}}$ and taking a diagonal subsequence 
of the double sequence $(i, j)$ gives a maximal limit orbifold domain 
$(\Omega, g)$, having the structure described in Theorem 4.1. As noted above, 
$\Omega$ may be compact, in which case $\Omega$ is denoted by $V$ and $V$ 
is an orbifold associated to $M$. At the other extreme, one may have 
$\Omega  = \emptyset$; in this case $M^{\nu_{0}} = \emptyset$, for any given 
$\nu_{0}$, provided $i$ is sufficiently large. 

 Next, consider the $\nu_{0}$-thin part of $(M, g_{i})$, i.e. 
\begin{equation} \label{e4.11}
M_{\nu_{0}} = \{x\in M: \nu (x) \leq  \nu_{0}\}. 
\end{equation}
By the discussion above, we may assume that $\nu_{0}$ is arbitrarily 
small, in that $\nu_{0} = \nu_{0}(i) \rightarrow 0$ sufficiently 
slowly, as $i \rightarrow  \infty$. We divide $M_{\nu_{0}}$ into two 
further subdomains. Thus, fix a large constant $K < \infty$, and let
\begin{equation} \label{e4.12}
U = \{x\in M_{\nu_{0}}:\rho (x) > K\nu (x)\}, \ {\rm and} \ W = \{x\in 
M_{\nu_{0}}:\rho (x) \leq  K\nu (x)\}. 
\end{equation}
Thus $W$ corresponds to the set where the curvature may concentrate in $L^{2}$, 
in the scale of the volume radius, cf. also [4, p.63] for the analogous 
description in the case of Einstein metrics. 
\begin{lemma} \label{l4.2}
In a subsequence, the set $W$ tends metrically to finitely many points 
$\{z_{j}\}$, as $i \rightarrow \infty$. 
\end{lemma}

{\bf Proof:}
For $z_{i} \in W$, rescale the metrics $g_{i}$ so that $\nu (z_{i}) = 1$, 
i.e. set $\hat g_{i} = \nu (z_{i})^{-2}g_{i}$, so that $\hat \rho(z_{i}) 
\leq K$. By the definition of $L^{2}$ curvature radius in 
(4.4), one has on $\hat g_{i}$, 
\begin{equation} \label{e4.13}
\int_{B(\hat \rho)}|\hat R|^{2} = 
c_{0}\frac{volB(\hat \rho)}{\hat \rho^{4}}, 
\end{equation}
where $\hat \rho = \hat \rho(z_{i})$, $B(\hat \rho) = B_{z_{i}}(\hat \rho)$. 
If $\hat \rho \leq 1$, then $vol B(\hat \rho) \geq \mu \hat \rho^{3}$, since 
$\hat \nu = 1$, while if $\hat \rho > 1$, $vol B(\hat \rho) \geq vol B(1) \geq \mu$. 
Thus
\begin{equation} \label{e4.14}
\int_{B(\hat \rho)}|\hat R|^{2} \geq  c_{1}, 
\end{equation}
with $c_{1}$ depending only on $c_{0}$. $\mu$, and $K$. Hence, by 
scale-invariance and (4.9), there is a bounded number of such balls, 
with bound depending only on $c_{0}$, $\mu$, $K$ and $\Lambda'$. 
Since $\nu_{i}(z_{i}) \rightarrow 0$, as $i \rightarrow \infty$, 
$\rho (z_{i}) \rightarrow 0$, and so these balls converge metrically 
to a finite number of points as $i \rightarrow \infty$. 
{\endproof}

 Consider now the complement $U$. For any given $x_{i}\in U$, we work 
again in the scale $\hat g_{i}$ where $\nu = 1$, so that 
$\hat \rho(x_{i}) \geq K$. In this scale, the bound (4.1) becomes
\begin{equation} \label{e4.15}
\int_{M}|\hat Ric|^{2p} \leq  \Lambda\nu_{i}^{4p-4} \rightarrow 0 \ \ 
{\rm as} \ \ i \rightarrow \infty.  
\end{equation}
When the $L^{2p}$ norm of the Ricci curvature is sufficiently small, 
one has a volume comparison result, cf. [36] for instance, which gives
\begin{equation} \label{e4.16}
\hat \nu (y_{i}) \geq  \mu_{1}\hat \nu (x_{i}) = \mu_{1}, 
\end{equation}
for all $y_{i}\in  B_{x_{i}}(K)$, where $\mu_{1}$ depends only on $K$ and 
$\mu$. Thus, the ball $B_{x_{i}}(K)$ in the metric $\hat g_{i}$ is 
everywhere non-collapsed as $i \rightarrow \infty$. By (4.7), the $L^{2,2p}$ 
harmonic radius is uniformly bounded below on $B_{x_{i}}(K)$, and hence the 
metrics $\hat g_{i}$ are precompact in the $L^{2,2p}$ and $C^{\alpha}$ 
topologies. As discussed following (4.7), $c_{0}$ is chosen sufficiently 
small so that the metric is close to the flat metric. Since 
$\hat \nu (x_{i})  = 1$ and $K$ is large, each $B_{x_{i}}(K)$ is thus 
close to a ball of radius $K$ in a non-trivial flat manifold 
${\mathbb R}^{4}/ \Gamma$, where $\Gamma$ is a discrete group of Euclidean 
isometries acting freely on ${\mathbb R}^{4}$. 

 Hence, associated to every point in $U$, there is a neighborhood of a 
definite size, which is diffeomorphic to ${\mathbb R}^{4}/\Gamma$, for 
some $\Gamma \neq \{e\}$, and on which one has uniform bounds on the 
metric in local harmonic coordinates in $C^{\alpha}\cap L^{2,2p}$. 
Following the collapse theory of Cheeger-Gromov [15], it is shown in [5] 
that these elementary $F$-structures piece together to give a global 
$F$-structure on $U$. (In [5], this is done in dimension 3, but the proof 
given works the same in all dimensions, given $C^{\alpha}\cap L^{2,q}$ 
control of the metric, for some $q > n = dim M$). 

  Finally, collapse at $\{x_{i}\}$ with locally bounded curvature in $L^{2p}$ 
means $\nu(x_{i}) \rightarrow 0$ and the scale-invariant quanitity 
$(\nu(x_{i}))^{4p-4}\int|R|^{2p} \rightarrow 0$, as $i \rightarrow \infty$. 
This follows from the results established above. 
{\endproof}

\begin{remark} \label{r 4.3}
  {\rm {\bf (i).} An essentially immediate consequence of the proof of Theorem 
4.1 is that if 
$$\inf_{{\mathbb M}_{1}}\int|Ric|^{2p} = 0,$$
then either $M$ admits a Ricci-flat orbifold metric, i.e. $\Omega = M$ or 
$\Omega = V$, or if not, then $\Omega = \emptyset$ so that $M$ carries an 
$F$-structure on the complement of finitely many, arbitrarily small metric balls. 
To see this, suppose $\Omega \neq \emptyset$. We then claim  $(\Omega, g)$ is 
necessarily complete. The incompleteness of $(\Omega, g)$ is caused (only) by 
the collapse of the metric in finite distance. However, $(\Omega, g)$ is Ricci-flat, 
and the Bishop-Gromov volume comparison theorem rules out collapse of the volumes 
of local balls within finite distance. This proves that $(\Omega, g)$ is complete. 
On the other hand, a result of Calabi and Yau implies that a complete Ricci-flat metric 
on an open manifold has infinite volume, contradicting (4.2).  

  {\bf (ii).}
 The lower semi-continuity of the functional implies of course that $\int_{\Omega}|Ric|^{2p} 
\leq \Lambda$. In local harmonic coordinates, the metric $g$ is $C^{\alpha}$ and satisfies 
(4.8). Hence
$$R \in L_{loc}^{2p},$$
on the regular set $\Omega_{0}$ of $\Omega$. 

  {\bf (iii).}
 We point out if the bound (4.1) is replaced by the stronger bound 
$\int|R|^{2p} \leq  \Lambda$, then the same proof shows that the limit 
$(\Omega, g)$ has no orbifold singularities, and there are no 
singularities in the collapsed part either. This is because all 
singularities under the bound (4.1) arise from rescalings (blow-downs) 
of complete, non-compact Ricci-flat 4-manifolds, (ALE in case of orbifold 
singularities). Under a bound on $\int|R|^{2p}$, these cannot arise, since all 
blow-up limits under such a bound are necessarily flat, cf. (4.15). }
\end{remark}

  Now consider perturbations of one of the functionals ${\mathcal F}$ in 
(3.1) in the direction of the $L^{2p}$ norm of the Ricci curvature, 
(for example). We first consider ${\mathcal F}  = {\mathcal R}ic^{2}$, and set 
\begin{equation} \label{e4.17}
{\mathcal R}ic_{\varepsilon}^{2p} = \int|Ric|^{2} + \varepsilon\int 
(1+|Ric|^{2})^{p}, 
\end{equation}
for $p > 1$ and $\varepsilon > 0$. Similar perturbations of ${\mathcal W}^{2}$ 
and ${\mathcal R}^{2}$ will be discussed later. The perturbation 
(4.17) is analogous to the $\alpha$-energy perturbation of 
Sacks-Uhlenbeck [38] in their study of harmonic maps of surfaces into 
Riemannian manifolds; perturbations of this type for curvature 
functionals were studied in detail in dimension 3 in [5], [7]. 

 The functional ${\mathcal R}ic_{\varepsilon}^{2p}$ is a $C^{\infty}$ smooth 
functional on space ${\mathbb M}_{1}$ of metrics of volume 1. The idea is 
to obtain an existence result for minimizers of ${\mathcal R}ic_{\varepsilon}^{2p}$, 
and then pass to a limit $\varepsilon \rightarrow 0$, (or $p \rightarrow 1$), 
to obtain an existence result for minimizers of ${\mathcal R}ic^{2}$. To begin, 
the following is an immediate consequence of Theorem 4.1:

\begin{corollary} \label{c 4.4.}
  For any $\varepsilon > 0$ and $p > 1$, there exists a minimizing pair 
$(\Omega_{\varepsilon}, g_{\varepsilon})$ for ${\mathcal R}ic_{\varepsilon}^{2p}$. 
The pair $(\Omega_{\varepsilon}, g_{\varepsilon})$ has the properties given in 
Theorem 4.1 together with 
\begin{equation} \label{e4.18}
{\mathcal R}ic_{\varepsilon}^{2p}(g_{\varepsilon}) \leq  
\inf_{{\mathbb M}_{1}}{\mathcal R}ic_{\varepsilon}^{2p}. 
\end{equation}
\end{corollary}
{\endproof}

 The Euler-Lagrange equations of ${\mathcal R}ic_{\varepsilon}^{2p}$ on ${\mathbb M}_{1}$ 
at a metric $g$ are:
\begin{equation} \label{e4.19}
D^{*}Dh - 2\delta^{*}\delta h - (\delta\delta h) g + {\mathcal P}_{R}  = 0, 
\end{equation}
where $h = fRic$, $f = (1+\varepsilon p(1+|Ric|^{2})^{p-1})$ and ${\mathcal P}_{R}$ 
is a curvature term, given by
\begin{equation} \label{e4.20}
{\mathcal P}_{R}  = - 2R(h) + \frac{1}{2}[|Ric|^{2}+\varepsilon 
(1+|Ric|^{2})^{p}+ c_{\varepsilon}]g, 
\end{equation}
where $c_{\varepsilon}$ is a constant; ($c_{\varepsilon} = 
\varepsilon \inf_{{\mathbb M}_{1}}[(p-1)(\int (1+|Ric|^{2}))^{p} + 
p(\int (1+|Ric|^{2}))^{p-1})$). We note that $(\delta^{*}\omega)(X,Y) = 
\frac{1}{2}[(\nabla_{X}\omega)(Y) + (\nabla_{Y}\omega)(X)]$, so that 
$\delta \omega = -tr \delta^{*}\omega$. 

 This formula is easily derived from standard formulas, cf. [12]. Briefly, 
the variation of the integrand of ${\mathcal R}ic_{\varepsilon}^{2p}$ is 
$$(1+2\varepsilon p(1+|Ric|^{2})^{p-1})[2\langle Ric'(g'), Ric \rangle - 
2\langle Ric^{2}, g' \rangle] + \tfrac{1}{2}[|Ric|^{2} + \varepsilon(1+|Ric|^{2})^{p} + c] 
\langle g', g \rangle.$$
One has
\begin{equation} \label{e4.21}
2Ric'(k) = D^{*}Dk - 2\delta^{*}\delta k - D^{2}(tr k) - 2R(k) + Ric \circ k + 
k \circ Ric, 
\end{equation}
and if $(Ric')^{*}$ denotes the adjoint of $Ric'$, then 
\begin{equation} \label{e4.22}
2(Ric')^{*}(k) = D^{*}Dk - 2\delta^{*}\delta k - (\delta\delta k)g - 2R(k) 
+ Ric \circ k + k \circ Ric.
\end{equation}
Combining the various terms gives (4.19). 

  When $\varepsilon = 0$, $h = Ric$ and (4.19) 
becomes the Euler-Lagrange equation $\nabla{\mathcal R}ic^{2} = 0$, i.e. 
\begin{equation} \label{e4.23}
D^{*}DRic - 2\delta^{*}\delta Ric - (\delta\delta Ric)g - 2R(h) + 
\tfrac{1}{2}|Ric|^{2}g = 0. 
\end{equation}

 Now a minimizer $(\Omega_{\varepsilon}, g_{\varepsilon})$ from Corollary 4.4 
is a weak $L^{2,2p}$ solution of the Euler-Lagrange equation (4.19), 
i.e. (4.19) holds when viewed as a distribution and paired with any 
$L^{2,2p^{*}}$ symmetric bilinear form $\psi$, of compact support in 
$\Omega_{\varepsilon}$; $2p^{*} = 1 - \frac{1}{2p}$ is the conjugate 
exponent to $2p$. Note that $h \in L^{2p/2p-1}$, since $Ric \in  L^{2p}$. 
For the remainder of this section, we assume $p > 1$ with $p$ close to 1. 

  In fact, the regularity can be improved a little for metrics $g$ which are local 
minima of ${\mathcal R}ic_{\varepsilon}^{2p}$, in the sense that $g$ minimizes 
${\mathcal R}ic_{\varepsilon}^{2p}$ among nearby, compact perturbations of $g$. 
 
\begin{lemma} \label{l 4.5}
Let $g$ be an $L^{2,2p}$ metric which locally minimizes ${\mathcal R}ic_{\varepsilon}^{2p}$ 
on a domain $U$, with $\varepsilon$ small. Then, 
\begin{equation} \label{e4.24}
Ric \in L_{loc}^{4},
\end{equation}
so that $g$ is locally in $L^{2,4}$. 
\end{lemma}
{\bf Proof:} This is proved in [46] for metrics locally minimizing the $L^{2p}$ norm 
of the curvature $R$, and the proof for ${\mathcal R}ic_{\varepsilon}^{2p}$ is essentially 
the same; thus we will be somewhat brief and refer to [46] for further details. 

  The local Ricci flow $\frac{d}{dt}g(t) = -2\chi g(t)$, where $\chi$ is a local 
cutoff function supported in $U$, is defined for metrics in $L^{2,2p}$. Using the 
local minimizing property of $g = g(0)$, it suffices to show that 
$$\frac{d}{dt}{\mathcal R}ic_{\varepsilon}^{2p}(g(t))|_{t=0} + c(\int\chi^{2}|Ric|^{4})^{1/2} 
\leq C[1 + {\mathcal R}ic_{\varepsilon}^{2p}(g)].$$
It follows from (4.21) and the Bianchi identity $\delta Ric = -\frac{1}{2}ds$, 
that 
$$\frac{d}{dt}{\mathcal R}ic_{\varepsilon}^{2p}(g(t)) \leq 
-\int \langle D^{*}D\chi Ric, fRic \rangle + C[1+\int \chi |Ric|^{2p}|R|],$$ 
where $f$ is defined as following (4.19). It is important for the following to 
observe that $f \geq 1$. One has $-\int \langle D^{*}D\chi Ric, fRic \rangle = 
-\int \langle D\chi Ric, DfRic \rangle$. Expanding this out and using the 
Cauchy-Schwarz and Young inequalities, ($ab \leq \mu a^{2} + \mu^{-1}b^{2}$), 
together with the fact that $\varepsilon$ is small, gives  
$-\int \langle D\chi Ric, DfRic \rangle \leq - \int \chi |DRic|^{2} + 
C[1+\int \chi |Ric|^{2p+1}]$, so that  
$$\frac{d}{dt}{\mathcal R}ic_{\varepsilon}^{2p}(g(t))_{t=0} + 
c\int \chi|DRic|^{2} \leq  C[1 + \int \chi |Ric|^{2p}|R|].$$
Since $|d|Ric||^{2} \leq c|DRic|^{2}$, the result then follows easily from the 
local Sobolev embedding, $L^{4} \subset L^{1,2}$ in dimension 4, via use of 
Remark 4.3(ii) and the H\"older inequality on the term 
$|Ric|^{2p}|R| \leq |Ric|^{2}|R|^{2p-1}$. 
{\endproof}

 An important point at this stage is to see that weak solutions of 
(4.19) are smooth. 

\begin{proposition} \label{p 4.6.}
  Any locally defined weak $L^{2,2p}$ solution $g$ of (4.19), which is a local 
minimizer of ${\mathcal R}ic_{\varepsilon}^{2p}$ is $C^{\infty}$ smooth. 
\end{proposition}
{\bf Proof:}
  The idea is of course to use elliptic regularity results to boost the regularity 
of $g$. By (4.24) and Remark 4.3(ii), the term ${\mathcal P}_{R} \in L^{4/2p}$, so that 
(4.19) has the form 
\begin{equation} \label{e4.25}
L(h) = D^{*}Dh - 2\delta^{*}\delta h - (\delta \delta h)g = -{\mathcal P}_{R} \in L^{4/2p}.
\end{equation}
However, the linear operator $L$ is not elliptic, due to its invariance properties 
under the diffeomorphism group. To deal with this, any symmetric bilinear form 
$\psi \in T_{g}{\mathbb M}$ on $M$ may be decomposed as 
\begin{equation} \label{e4.26}
\psi = \delta^{*}X + \phi g + k,
\end{equation}
where $X$ is a vector field, $\phi$ a function and $k$ is transverse-traceless, 
$\delta k = tr k = 0$. This also holds locally, and so $h = fRic \in L^{p'}$, 
$p' = 4/(2p-1)$ may be written as a sum as in (4.26), with $\delta^{*}X, \phi$ 
and $k$ in $L^{p'}$.

  The equation (4.25) thus decomposes as 
\begin{equation} \label{e4.27}
L(k) + L(\delta^{*}X) + L(\phi g) \in L^{4/2p}.
\end{equation}
One has $L(k) = D^{*}Dk$, which is of course an elliptic 
operator on $k$. To compute $L(\delta^{*}X)$, note that ${\mathcal L}_{X}(Ric) = 
Ric'(\delta^{*}(X))$, while in local coordinates, ${\mathcal L}_{X}(Ric) \sim 
X(Ric) + (\partial X)Ric$. The second term is in $L^{4/2p}$ while the first term 
is in $L^{-1,q}$, where $q = 4$. Thus ${\mathcal L}_{X}(Ric) \in L^{-1,q}$. 
On the other hand, one has $2(Ric')^{*}(\beta) = 2Ric'(\beta) + D^{2} tr \beta - 
(\delta \delta \beta) g$, (see (4.21)-(4.22)). It follows that 
$$L(\delta^{*}X) = -\delta \delta \delta^{*}X g - D^{2}\delta X + L^{-1, q},$$
where have used $tr \delta^{*}X = -\delta X$, and $+ L^{-1, q}$ denotes the 
addition of an element in $L^{-1,q}$. A simple direct computation gives  
$$L(\phi g) = -2\Delta \phi g + 2D^{2}\phi.$$
Now by a standard Weitzenbock formula, 
$$\delta \delta^{*}X = D^{*}D X = \delta d X + d \delta X - Ric(X),$$
so that $\delta \delta \delta^{*}X = \delta d \delta X - \delta(Ric(X))$. 
One may decompose $X$ as 
$$X = \nabla \omega + Y,$$
where $\delta Y = 0$. In sum, combining these computations gives
\begin{equation} \label{e4.28}
L(\delta^{*}X) + L(\phi g) = -(\Delta \Delta \omega + 2\Delta \phi)g + 
D^{2}(\Delta \omega + 2\phi) + L^{-1, q}.
\end{equation}

  Now consider first the trace of the equation (4.27). Since $k$ is 
transverse-traceless, via (4.28) this gives
$$\Delta \Delta \omega + 2\Delta \phi \in L^{-1, q}.$$
The coefficients of the Laplacian $\Delta$ are in $C^{\alpha}$ in local harmonic 
coordinates. By elliptic regularity, cf. [31], it follows that $\Delta \omega + 
2\phi \in L^{1,q}$, and hence 
$$L(\delta^{*}X) + L(\phi g) \in L^{-1, q}.$$
Returning to (4.27), we then have $L(k) \in L^{-1,q}$, so that again by 
elliptic regularity, 
\begin{equation} \label{e4.29}
k \in L^{1,q},
\end{equation}
giving the main initial regularity boost. 

  Since $\Delta \omega, \phi \in L^{p'}$, while $\Delta \omega + 2\phi \in L^{1,q}$, 
if $\Delta \omega$ and $\phi$ are linearly independent, it follows that $\Delta \omega 
\in L^{1,q}$ and $\phi \in L^{1,q}$. In fact, this is case unless $\Delta \omega + 
2\phi \equiv 0$. Suppose first that $\Delta \omega + 2\phi \neq 0$, so that 
$\Delta \omega \in L^{1,q}$, $\phi \in L^{1,q}$, and hence $tr h = \Delta \omega 
+ 4\phi \in L^{1,q}$. Thus, one has
\begin{equation} \label{e4.30}
h = \delta^{*}Y + L^{1,q}.
\end{equation}
To show that $h = fRic \in L^{1,q}$, take the exterior derivative of (4.30), giving
$$dh = f dRic + df \wedge Ric = d\delta^{*}Y + L^{q} = R(Y) + L^{q} \in L^{q'},$$
where the terms are 2-forms with values in the tangent bundle and $q' < 4$. 
Taking the trace on the last two indices and using the Bianchi identity 
$\delta Ric = -\frac{1}{2}ds$, one obtains
$$tr dh= d tr h + \delta h = \frac{1}{2}d(tr h) - E(df) \in L^{q'},$$
and so $\delta h \in L^{q'}$, since $d tr h \in L^{q}$. 

  On the other hand, from (4.26) with $\psi = h$ and the Weitzenbock formula, one has 
$\delta h = \delta d Y - d\Delta \omega - d\phi - Ric(X) = \delta dY + L^{q'}$. 
It follows that $\delta d Y \in L^{q'}$. Since $\delta Y = 0$, this gives $Y \in L^{2,q'}$, 
and hence $h \in L^{1,q'}$. In turn, this now implies $Ric \in L^{1,(2p-1)q'}$, so that 
$g \in L^{3,(2p-1)q'} \subset C^{2,\alpha}$ in local harmonic coordinates. 

  Suppose instead $\Delta \omega + 2\phi = 0$. Then in place of (4.30) one has 
$$h = \delta^{*}X - \tfrac{1}{2}(\Delta \omega)g + L^{1,q}.$$
Via the Weitzenbock formula, this gives $\delta \delta h = -\Delta \Delta \omega + 
L^{-1,q}$. Also, as below (4.30), one has $d tr h + \delta h = -\frac{3}{2}d\Delta \omega 
+ L^{q'}$, so that $-\Delta tr h + \delta \delta h = \frac{3}{2}\Delta \Delta \omega 
+ L^{-1, q'}$. Since $tr h = -\Delta \omega$, this gives $\Delta \Delta \omega \in 
L^{-1,q'}$, and hence $\Delta \omega \in L^{q'}$. It follows then as before that 
$h \in L^{1,q'}$. 

  One may now repeat the arguments above inductively, improving the regularity of $h$ and 
$g$ at each step, to obtain $g \in C^{\infty}$. 
{\endproof}

\begin{conjecture} \label{conj4.7}
 {\rm The minimizers of ${\mathcal R}ic_{\varepsilon}^{2p}$ are complete, i.e. 
$g_{\varepsilon}$ is complete on each component of $\Omega_{\varepsilon}$, and 
$\Omega_{\varepsilon}$ has only finitely many components. 

 Results of this type are proved for minimizers of analogous 
functionals in dimension 3 in [5]. However, it seems difficult to 
extend the proof in the 3-dimensional case to 4-dimensions. }
\end{conjecture}

\medskip

 Consider now the analogous procedure for either ${\mathcal R}^{2}$ or 
${\mathcal W}^{2}$, i.e. the perturbations
\begin{equation} \label{e4.31}
\int|R|^{2} + \varepsilon\int (1+|R|^{2})^{p}, \ \ {\rm or} \ \  \int|W|^{2} + 
\varepsilon\int (1+|W|^{2})^{p}. 
\end{equation}
Corollary 4.4 obviously holds for the first of these perturbed 
functionals, (via Theorem 4.1), but it is unknown if it holds for the 
second. Even if it did, it is not clear if Proposition 4.6 holds for either 
of the functionals in (4.31). For example, the Euler-Lagrange equation for the 
perturbation of ${\mathcal R}^{2}$ has a term of the form
$$(1+\varepsilon p(1+|R|^{2})^{p-1}) \langle R', R \rangle = 
(R')^{*}[(1+\varepsilon p(1+|R|^{2})^{p-1}R]. $$
This has a much more complicated form than (4.19), and the proof of 
Proposition 4.6 will not apply directly. The same remarks apply to the 
functional ${\mathcal W}^{2}$. A similar difficulty remains if one perturbs 
either functional by $(1+|Ric|^{2})^{p}$ in place of $(1+|R|^{2})^{p}$ 
or $(1+|W|^{2})^{p}$. Instead, we use a slightly different path. 

\medskip

 By the Chern-Gauss-Bonnet theorem (3.3), minimizing ${\mathcal R}^{2}$ is the 
same as minimizing ${\mathcal Z}^{2}$; these two functionals also have the 
same Euler-Lagrange equations or critical points. Similarly, minimizers 
or critical points of ${\mathcal W}^{2}$ are the same as those of ${\mathcal Z}^{2}$ - 
$\frac{1}{12}{\mathcal S}^{2} = {\mathcal R}ic^{2}$ - $\frac{1}{3}{\mathcal S}^{2}$, 
cf. (3.4). Thus, consider first the functional
\begin{equation} \label{e4.32}
{\mathcal Z}_{\varepsilon}^{2p}(g) = \int(|Ric|^{2} - \frac{1}{4}s^{2}) + 
\varepsilon\int (1+|Ric|^{2})^{p}. 
\end{equation}
For $\varepsilon > 0$ small, this is a small perturbation of ${\mathcal Z}^{2}$, with the 
property that a bound on ${\mathcal Z}_{\varepsilon}^{2p}$ implies a bound on 
${\mathcal R}ic^{2p}$, as in (4.1). 

  The Euler-Lagrange equation for ${\mathcal Z}_{\varepsilon}^{2p}$ is 
\begin{equation} \label{e4.33}
D^{*}Dh - 2\delta^{*}\delta h - (\delta\delta h)g - \tfrac{1}{2}
(D^{2}s - (\Delta s) g) + {\mathcal P}_{Z} = 0, 
\end{equation}
where $h = fRic$, $f = (1+\varepsilon p(1+|Ric|^{2})^{p-1})$, 
with ${\mathcal P}_{Z}$ given by
\begin{equation} \label{e4.34}
{\mathcal P}_{Z} = - 2R(h) + \tfrac{1}{2}sRic + \frac{1}{2}(\zeta + c_{\varepsilon})g,
\end{equation}
where $\zeta = |Ric|^{2} + \varepsilon (1+|Ric|^{2})^{p}$. 

  To derive (4.33), for $f$ and $\zeta$ as above, the variation of the integrand of 
${\mathcal Z}_{\varepsilon}^{2p}$ is 
$$2f[\langle Ric', Ric \rangle  - \langle Ric^{2}, g' \rangle] - \tfrac{1}{2}ss' 
+ \tfrac{1}{2}(\zeta + c_{\varepsilon}) \langle g', g \rangle .$$
One has $2f\langle Ric', Ric \rangle = 2(Ric')^{*}(fRic)$, and similarly $ss' = (s')^{*}(s)$. 
The first term is given by (4.22) while the second term is $D^{2}s - (\Delta s)g - s Ric$, 
and combining these expressions gives (4.33). 

 As before, when $\varepsilon = 0$, $h = Ric$ and (4.33) becomes the 
Euler-Lagrange equation $\nabla{\mathcal Z}^{2} = 0$, i.e. 
\begin{equation} \label{e4.35}
D^{*}Dz - 2\delta^{*}\delta z - (\delta\delta z) g - 2R(z) + \tfrac{1}{2}|z|^{2}g = 0. 
\end{equation}

\begin{corollary} \label{c 4.8.}
   Corollary 4.4 and Proposition 4.6 hold for the functional ${\mathcal Z}_{\varepsilon}^{2p}$. 
\end{corollary}
{\bf Proof:}
 The proof of Corollary 4.4 via Theorem 4.1 is the same as before, using the fact noted 
above that a bound on ${\mathcal Z}_{\varepsilon}^{2p}$ implies a bound on $Ric^{2p}$. 
It is straightforward to show that Lemma 4.5 also holds for ${\mathcal Z}_{\varepsilon}^{2p}$ 
in place of ${\mathcal R}ic_{\varepsilon}^{2p}$. The only difference is that the term 
$|DRic|^{2}$ is replaced by $|DRic|^{2} - \frac{1}{4}|ds|^{2}$. However, $Ric = z + 
\frac{s}{4}g$ and $|Dz|^{2} \geq \frac{1}{4}(\delta z)^{2} = \frac{1}{64}|ds|^{2}$, so 
that $|DRic|^{2} - \frac{1}{4}|ds|^{2} \geq \frac{1}{64}|ds|^{2}$. Given this adjustment, 
the rest of the proof carries over as before. 

  Moreover, the proof of Proposition 4.6 is also essentially the same as before. The 
addition of the $s$ terms in (4.33) implies that (4.28) is changed to 
\begin{eqnarray} \label{e4.36}
L(\delta^{*}X) + L(\phi g) - \tfrac{1}{2}(D^{2} s - \Delta s g)  \\
= D^{2}(\Delta \omega + 2\phi -\tfrac{1}{2}s) - 
(\Delta \Delta \omega + 2\Delta \phi - \tfrac{1}{2}\Delta s)g . \nonumber
\end{eqnarray}
Given this, the rest of the proof remains valid, with only minor changes. 
{\endproof}

 Next, we treat the case of ${\mathcal W}^{2}$ in a similar way. Among various possible 
perturbations, consider 
\begin{equation} \label{e4.37}
{\mathcal W}_{\varepsilon,\lambda}^{2p}(g) = \int(|Ric|^{2} - \frac{1}{3}s^{2}) + 
\varepsilon\int (1+ |Ric|^{2})^{p} + \lambda \int s^{2} = {\mathcal R}ic_{\varepsilon}^{2p} 
- (\frac{1}{3} - \lambda){\mathcal S}^{2},
\end{equation}
which, modulo $\chi (M)$, gives $2{\mathcal W}^{2}$ when $\varepsilon = \lambda = 0$. 
As will be seen later, the relations between $\varepsilon > 0$, $p > 1$ and $\lambda > 0$ 
determine the choice of conformal gauge as $\varepsilon \rightarrow 0$. It is possible 
in the following to dispense with $\lambda$, i.e. set $\lambda = 0$, but we do not do so 
since a suitable choice of $\lambda > 0$ leads to a Yamabe-type gauge in the limit 
$\varepsilon = 0$. As discussed in Example 3.1, a good choice of conformal gauge is 
of some importance. It is clear that a bound on ${\mathcal W}_{\varepsilon,\lambda}^{2p}$ 
implies a bound on ${\mathcal R}ic^{2p}$. 

  The Euler-Lagrange equation for ${\mathcal W}_{\varepsilon,\lambda}^{2p}$ is 
\begin{equation} \label{e4.38}
D^{*}Dh - 2\delta^{*}\delta h - (\delta\delta h)g - 2(\tfrac{1}{3} - \lambda)
(D^{2}s - (\Delta s) g) + {\mathcal P}_{W} = 0,
\end{equation}
where $h = fRic$ is as in (4.19). The curvature term ${\mathcal P}_{W}$ is 
\begin{equation} \label{e4.39}
{\mathcal P}_{W}= -2R(h) + 2(\tfrac{1}{3} - \lambda)sz + \tfrac{1}{2}[|Ric|^{2} + 
\varepsilon(1+|Ric|^{2})^{p} + c_{\varepsilon}]g. 
\end{equation}
The derivation of (4.38) is straightforward, given the derivations of (4.19) and (4.33). 
A little computation shows that the trace of (4.38) is given by
\begin{equation} \label{e4.40}
-\Delta trh - 2\delta\delta h + (2 - 6\lambda)\Delta s - 2\varepsilon p
(1 + |Ric|^{2})^{p-1}|Ric|^{2} + 2[\varepsilon(1+|Ric|^{2})^{p} + c_{\varepsilon}].
\end{equation}

   When $\varepsilon = \lambda = 0$, $h = Ric$ and $f = 1$, and so (4.38) becomes
\begin{equation} \label{e4.41}
D^{*}DRic + \tfrac{1}{3}D^{2}s + \tfrac{1}{6}\Delta s g - 2R(Ric) + \tfrac{2}{3}sRic + 
\tfrac{1}{2}[|Ric|^{2} - \tfrac{1}{3}s^{2}]g = 0. 
\end{equation}
The equation (4.41) is of course the conformally invariant Bach equation, i.e. 
the Euler-Lagrange equations $\nabla{\mathcal W}^{2} = 0$ for ${\mathcal W}^{2}$, cf. [12]. 

 We are now in position to prove:
\begin{corollary} \label{c 4.9.}
   Corollary 4.4 and Proposition 4.6 hold for the functional ${\mathcal W}_{\varepsilon,\lambda}^{2p}$.  
\end{corollary}
{\bf Proof:}
 As noted above, a bound on ${\mathcal W}_{\varepsilon,\lambda}^{2p}$ implies a bound 
on ${\mathcal R}ic^{2p}$, so Theorem 4.1 and Corollary 4.4 hold for 
${\mathcal W}_{\varepsilon,\lambda}^{2p}$.

  To verify that (4.24) holds, write ${\mathcal W}_{\varepsilon,\lambda}^{2p}$ as 
$${\mathcal W}_{\varepsilon,\lambda}^{2p} = 2{\mathcal W}^{2} - 16\pi^{2}\chi(M) + 
\varepsilon \int (1 + |Ric|^{2})^{p} + \lambda \int s^{2}.$$
The three leading order terms in the Euler-Lagrange equation (4.41) for ${\mathcal W}^{2}$ 
may be written in the form $\delta d(Ric - \frac{s}{6}g)$. In place of the local Ricci flow, 
deform the metric locally in the direction $-2(Ric - \frac{s}{6}g)$, i.e. 
$\frac{d}{dt}g = -2\chi(Ric - \frac{s}{6}g)$. Using the results in [18], the proof 
that this local flow exists is the same as that for the local Ricci flow. As in the 
proof of Lemma 4.5, one then has 
$$\frac{d}{dt}{\mathcal W}^{2}(g(t)) \leq C(1 + \int \chi |R|^{3}) .$$
On the other hand, the same estimates as before in the proof of Lemma 4.5 hold for the 
variation of functional $\varepsilon \int (1 + |Ric|^{2})^{p}$ in the direction of the 
modified local Ricci flow, ($-2\chi(Ric - \frac{s}{6}g)$), in place of the variation of 
${\mathcal R}ic_{\varepsilon}^{2p}$ or ${\mathcal Z}_{\varepsilon}^{2p}$ along the local 
Ricci flow; (this uses the fact that $\frac{1}{6} < \frac{1}{4}$ from the proof of 
Corollary 4.8). The variation of $\lambda {\mathcal S}^{2}$ contributes only lower order 
terms, and so the estimates are unaffected. One then obtains (4.24) in the same way as 
before. 

  It is also straightforward to see that the proof of Proposition 4.6 for 
${\mathcal W}_{\varepsilon,\lambda}^{2p}$ proceeds exactly as before in the case of 
${\mathcal R}ic_{\varepsilon}^{2p}$ or ${\mathcal Z}_{\varepsilon}^{2p}$. 
{\endproof}

  The same results hold for the functionals ${\mathcal W}_{\pm}^{2}$, which 
differ from ${\mathcal W}^{2}$ just by a topological term, (as with ${\mathcal R}^{2}$ and 
${\mathcal Z}^{2}$). 

\medskip

  Summarizing the work above, we have now produced suitable perturbations of the functionals 
${\mathcal F}$ in (3.1), and proved the existence of minimizing configurations for each of 
them. As mentioned above, the idea now is to take a sequence $\varepsilon  = 
\varepsilon_{i} \rightarrow 0$, and consider the behavior of 
(subsequences of) a sequence $(\Omega_{\varepsilon}, g_{\varepsilon})$ of 
minimizing pairs in the limit $\varepsilon \rightarrow 0$. Although the 
metrics $g_{\varepsilon}$ are $C^{\infty}$ smooth, one no longer has 
uniform control of the $L^{2p}$ norm of the Ricci curvature of 
$g_{\varepsilon}$, so that Theorem 4.1 does not apply. Nevertheless, using 
the Chern-Gauss-Bonnet theorem, together with the smoothness of 
$g_{\varepsilon}$ and the fact that $g_{\varepsilon}$ satisfies an (essentially) 
elliptic Euler-Lagrange equation, we show that the conclusions of Theorem 4.1 do 
in fact hold. As seen above and especially in Example 3.1, the cases 
${\mathcal R}ic^{2}$, ${\mathcal R}^{2}$ and the cases ${\mathcal W}^{2}$, ${\mathcal W}_{\pm}^{2}$ 
are somewhat different, and so these two situations are treated separately.

\begin{theorem} \label{t 4.10} (Geometric Decomposition with respect to 
${\mathcal R}ic^{2}, {\mathcal R}^{2}$.)
  
 Let $M$ be a closed, oriented 4-manifold and let ${\mathcal F}$ be one of 
the functionals ${\mathcal R}ic^{2}$, ${\mathcal R}ic^{2}$ on ${\mathbb M}_{1}$. 
Then minimizers of ${\mathcal F}$ on ${\mathbb M}_{1}$ are realized in the 
idealized sense that there exist minimizing sequences $\{g_{i}\}$ converging 
in the pointed Gromov-Hausdorff topology to one of the following configurations:

 {\rm (I)}. A compact, oriented, possibly reducible orbifold $(V, g_{0})$ associated to 
$M$, with $C^{\infty}$ metric $g_{0}$ on the regular set $V_{0}$. The metric $g_{0}$ 
extends $C^{0}$ across the orbifold singularities and 
\begin{equation} \label{e4.42}
vol_{g_{0}}V = 1. 
\end{equation}
Further
\begin{equation} \label{e4.43}
{\mathcal F} (g_{0}) +  \sum_{k}{\mathcal F}(N_{k}, g_{\infty}^{k}) = 
\inf_{g\in{\mathbb M}_{1}}{\mathcal F} (g),
\end{equation}
where the sum is over the collection of ALE spaces $(N_{k}, g_{\infty}^{k})$ 
associated with the singularities of $V$, cf. also (4.65). 

 {\rm (II)}. A maximal orbifold domain $\Omega \subset\subset M$, possibly 
reducible and possibly empty, with $C^{\infty}$ smooth metric $g_{0}$ on the 
regular set $\Omega_{0}$. The metric $g_{0}$ extends $C^{0}$ across the orbifold 
singularities, and satisfies
\begin{equation} \label{e4.44}
vol_{g_{0}}V \leq  1, 
\end{equation}
together with
\begin{equation} \label{e4.45}
{\mathcal F} (g_{0}) +  \sum_{k}{\mathcal F}(N_{k}, g_{\infty}^{k}) \leq  
\inf_{g\in{\mathbb M}_{1}}{\mathcal F} (g),
\end{equation}
where the sum is as in (4.43). 

   Outside a sufficiently large compact set $K \subset \Omega$, $\Omega$ 
carries an $F$-structure along which $g_{0}$ collapses with locally 
bounded curvature on approach to $\partial \Omega$. The complement 
$M \setminus K$ also carries an $F$-structure outside a finite number of balls 
$B_{z_{j}}(\varepsilon_{i})$ with respect to $g_{i}$, where $\varepsilon_{i} 
\rightarrow 0$ as $i \rightarrow \infty$. If $K_{i}$ is an 
exhaustion of $\Omega$, then $M\setminus K_{i}$ collapses 
everywhere, and collapses with locally bounded curvature away from the 
singular points $\{z_{j}\}$. 

 In all cases, the metric $g_{0}$ satisfies the Euler-Lagrange equation
\begin{equation} \label{e4.46}
\nabla{\mathcal F}  = 0, 
\end{equation}
and minimizes ${\mathcal F}$ among compact perturbations. 
\end{theorem}

  In comparing this result with Theorem 1.1, one should note that it is 
possible that $\Omega = \emptyset$. Thus, (II) above includes both cases 
(II) and (III) of Theorem 1.1. 

\medskip

{\bf Proof:} 
By Proposition 4.6 and Corollaries 4.4 and 4.8, there exist minimizing pairs 
$(\Omega_{\varepsilon}, g_{\varepsilon})$ for ${\mathcal F}$, for any given 
$\varepsilon > 0$. As in the proof of Theorem 4.1, consider a 
thick-thin decomposition of $\Omega_{\varepsilon}$. As before, the 
discussion below applies to a sequence $(\Omega_{\varepsilon_{i}}, 
g_{\varepsilon_{i}})$ with $\varepsilon_{i} \rightarrow 0$, but we will 
usually drop the subscript from the notation. For any given $\nu_{0} > 0$, 
let
\begin{equation} \label{e4.47}
\Omega_{\varepsilon}^{\nu_{0}} = \{x\in\Omega_{\varepsilon}: \nu (x) >  
\nu_{0}\}. 
\end{equation}
For any given $\rho_{0} > 0$, consider the subdomain 
\begin{equation} \label{e4.48}
\Omega_{\varepsilon}^{\nu_{0},\rho_{0}} = 
\{x\in\Omega_{\varepsilon}^{\nu_{0}}: \rho (x) >  \rho_{0}\}, 
\end{equation}
where $\rho$ is the $L^{2}$ curvature radius. 

 We claim that on $\Omega_{\varepsilon}^{\nu_{0},\rho_{0}}$, one has 
$C^{\infty}$ smooth convergence (in a subsequence) to a limit 
$\Omega_{0}^{\nu_{0},\rho_{0}}$, away from the boundary. The 
Euler-Lagrange equations for $g_{\varepsilon}$, (4.19) or (4.33), are 
(essentially) uniformly elliptic as $\varepsilon \rightarrow 0$. More 
precisely, although (4.19) or (4.25) is not elliptic, the proof of 
Proposition 4.6 shows that uniform elliptic estimates hold on each summand 
in (4.26), independent of $\varepsilon$. The same applies to the Euler-Lagrange 
equation (4.33) for ${\mathcal Z}_{\varepsilon}^{2p}$. 

  However, the elliptic regularity estimates require, to get started, control 
on the $C^{\alpha}$ norm of the leading order coefficients, for some fixed $\alpha > 0$. 
Thus, uniform estimates require uniform control of the metric locally in $C^{\alpha}$, 
(in harmonic coordinates). Since there is no longer a uniform bound on ${\mathcal R}ic^{2p}$, 
but only a uniform bound on ${\mathcal R}^{2}$, one thus needs some stronger initial control 
on $g_{\varepsilon}$ in order to proceed. We note that Lemma 4.5 is not uniform in 
$\varepsilon$, since, for instance, it requires uniform control on the local Sobolev 
constant. 

 To obtain this initial control, let $\rho^{(q)}$ be the $L^{q}$ curvature radius, 
$q > 2$, with fixed parameter $c_{0} = c_{0}(q)$ in (4.4). We first claim that if 
$c_{0} = c_{0}(2)$ is sufficiently small, then there is a constant 
$\delta_{0} > 0$, depending only on $c_{0}(2)$, $c_{0}(q)$ and 
$\nu_{0}$, such that 
\begin{equation} \label{e4.49}
\rho^{(q)}(x) \geq  \delta_{0}\rho (x), 
\end{equation}
for all $x \in \Omega_{\varepsilon}^{\nu_{0},\rho_{0}}$. This gives 
uniform local $L^{2,q}$ and so $C^{\alpha} \cap L^{1,k}$ control of the 
metric in harmonic coordinates on $\Omega_{\varepsilon}^{\nu_{0},\rho_{0}}$, 
with $\alpha  > 0$, $k > 4$. Thus, the proof of Proposition 4.6 applies 
uniformly as $\varepsilon \rightarrow 0$, which gives the claim above of 
smooth convergence. 

 The proof of (4.49) is by contradiction. If (4.49) is not true, then 
there exists $x \in   \Omega_{\varepsilon}^{\nu_{0},\rho_{0}}$ such that 
$\rho^{(q)}(x) \leq  \delta\rho (x),$ where $\delta $ is arbitrarily 
small. Without loss of generality, we may assume that $x$ realizes, or almost 
realizes, the minimal value of the ratio $\rho(x) / 
dist(x, \partial \Omega_{\varepsilon}^{\nu_{0},\rho_{0}})$. Work in the 
scale $\hat g_{\varepsilon} = \rho^{(q)}(x)^{-2}g_{\varepsilon}$ where 
$\hat \rho^{(q)}(x) = 1$, so that $\hat \rho(x) >> 1$ and $\nu (x) >> 1$. 
Now in this scale, one does have uniform local control of 
$\hat g_{\varepsilon}$ in $L^{2,q}$, (independent of $\varepsilon$), and 
so in $C^{\alpha}$ in harmonic coordinates. Hence, as noted above following 
(4.48), $\hat g_{\varepsilon}$ is thus uniformly controlled in $C^{k}$, 
for any $k$. 

 Now since $\hat \rho^{(q)}(x) = 1$, the metric $\hat g_{\varepsilon}$ has 
a definite amount of curvature in $L^{q}$ on $B_{x}(1)$, depending only on 
the choice of $c_{0}(q)$. However, $\rho (x) >> 1$, and $c_{0}(2)$ is very 
small, so that the curvature is very small in $L^{2}$ on $B_{x}(1)$. Since 
the metric $\hat g_{\varepsilon}$ is uniformly controlled in $C^{k}$, for any 
$k$, this is impossible, and so proves (4.49). Clearly, if $c_{0}(q)$ is fixed, 
the size of $c_{0}(2)$ may be explicitly estimated in terms of $c_{0}(q)$. 

 Next consider the complementary region of 
$\Omega_{\varepsilon}^{\nu_{0}}$ where $\rho \leq \rho_{0} \leq \nu_{0}$, where 
$\rho_{0}$ is (arbitrarily) small, i.e. 
\begin{equation} \label{e4.50}
(\Omega_{\varepsilon})^{\nu_{0}}_{\rho_{0}} = 
\{x\in\Omega_{\varepsilon}^{\nu_{0}}: \rho (x) \leq  \rho_{0}\}. 
\end{equation}
Then Lemma 4.2, in the scale $g_{\varepsilon}$, (not $\hat g_{\varepsilon}$), 
shows that there is at most a bounded number of points 
$q_{j} \in  \Omega_{\varepsilon}^{\nu_{0}}$ such that 
$(\Omega_{\varepsilon})^{\nu_{0}}_{\rho_{0}} \subset  B_{q_{j}}(\rho_{0}).$ 
Thus, there is bounded number of points, independent of $\varepsilon$, 
where the curvature $R$ can concentrate in $L^{2}$. In fact, the bound is 
independent of $\nu_{0}$, and depends only on $c_{0}$ and $\mu$. 

 It follows that if $\rho_{\varepsilon}$ is any sequence such that 
$\rho_{\varepsilon} \rightarrow 0$, as $\varepsilon \rightarrow 0$, 
$(\Omega_{\varepsilon})^{\nu_{0}}_{\rho_{\varepsilon}}$ converges metrically to 
a finite number of points. We claim that in the limit, each of these point 
singularities is an orbifold singularity, possibly reducible, and possibly with 
$\Gamma = \{e\}$. To prove this, we use the orbifold compactness theorem of [8], 
which, given the $L^{2}$ bound on the curvature $R$ from the Chern-Gauss-Bonnet 
theorem, and the lower bound $\nu_{0}$ on the volume radius, shows that 
it suffices to have the following small curvature estimate on 
$(\Omega_{\varepsilon}^{\nu_{0}}, g_{\varepsilon})$: if $\int_{B(r)}|R|^{2} 
\leq  \delta$, for some fixed $\delta$ small, then there is a constant $C$, 
independent of $\delta$, such that 
\begin{equation} \label{e4.51}
\sup_{B(\frac{1}{2}r)}|R|^{2} \leq  \frac{C}{volB(r)}\int_{B(r)}|R|^{2}. 
\end{equation}
Since we are working in $\Omega^{\nu_{0}}$ for a fixed (but arbitrary) 
$\nu_{0}$, this is equivalent to showing that the curvature is bounded 
in $L^{\infty}$ in balls of radius $\rho$, when such $\rho$ balls are scaled 
to radius $\rho = 1$; (the estimate in (4.51) is scale-invariant). However, this 
has already been done; (4.49) gives a lower bound on the $L^{q}$ curvature radius, 
and higher order control follows as before via elliptic regularity as in Proposition 
4.6. 

  Suppose there exists $\nu_{0} > 0$ such that $\nu(x) \geq \nu_{0}$, for all 
$x \in \Omega_{\varepsilon}$ as $\varepsilon = \varepsilon_{i} \rightarrow 0$. 
Then $\Omega_{\varepsilon}$ is a compact orbifold $V_{\varepsilon}$ associated 
to $M$ and the analysis above shows that $(V_{\varepsilon}, g_{\varepsilon})$ 
converges, in a subsequence, to a limit orbifold $(V, g_{0})$, smoothly away from 
orbifold singularities. The equation (4.42) is obvious from the smooth convergence 
while the equation (4.43) follows from the proof of the orbifold compactness theorem 
[8] above, together with the work in [9] and [11]; this issue is discussed further 
below, (following Theorem 4.15). This completes the proof in Case (I). 

  Next, suppose there exist $x_{i} \in \Omega_{\varepsilon}$ such that $\nu(x_{i}) 
\rightarrow 0$ as $\varepsilon \rightarrow 0$. One may then choose a sequence 
$\nu_{j} \rightarrow 0$ and consider the domain $(\Omega_{\varepsilon}^{\nu_{j}}, 
g_{\varepsilon})$. As in the proof of Theorem 4.1, a diagonal subsequence of $(i, j)$ 
converges in the pointed Gromov-Hausdorff topology to a limit maximal orbifold domain 
$(\Omega, g_{0})$. For the same reasons as above, the limit satisfies the properties 
claimed in Theorem 4.10. The equalities in (4.42) and (4.43) are replaced however 
by inequalities, since part of the volume and part of the value of 
${\mathcal F}(g_{\varepsilon})$ may be contained in the collapsing region. 

 To complete the proof, consider the structure of the complementary, collapsing region 
\begin{equation} \label{e4.52}
(\Omega_{\varepsilon})_{\nu_{0}} = \{x\in\Omega_{\varepsilon}: \nu (x) \leq  
\nu_{0}\}, 
\end{equation}
where $\nu_{0}$ is small, and may be assumed to be arbitrarily small 
for $i$ sufficiently large. As in (4.12), form the two subdomains
\begin{equation} \label{e4.53}
U_{\varepsilon} = \{x\in (\Omega_{\varepsilon})_{\nu_{0}}:\rho (x) > K\nu 
(x)\}, \ {\rm and} \  W_{\varepsilon} = \{x\in (\Omega_{\varepsilon})_{\nu_{0}}:\rho (x) 
\leq  K\nu (x)\}. 
\end{equation}
The proof that $W$ tends to finitely many points in Lemma 4.2 holds 
without any changes here also. Regarding $U_{\varepsilon}$, although (4.15) 
does not hold in the present situation, in the scale 
$\hat g_{\varepsilon} = \nu (x_{i})^{-2}g_{\varepsilon}$ where 
$\hat \nu(x_{i}) = 1$, one has $\hat \rho(x_{i}) \geq K$. 
Hence, as described above in and following (4.49), the metric $\hat g_{\varepsilon}$ 
is smoothly close to a flat metric on the ball $B_{x_{i}}(K-1)$. For the 
same reasons as in the proof of Theorem 4.1 following (4.16), such a ball 
is close to a ball in a non-trivial flat manifold ${\mathbb R}^{4}/\Gamma$, 
giving rise to an elementary $F$-structure. Again, as previously, these elementary 
$F$-structures patch together to give a global $F$-structure on $U_{\varepsilon}$. 
The same arguments based on Theorem 4.1 also show that there exist minimizing sequences 
$\{g_{i}\}$ for ${\mathcal F}$ on $M$ such that $M \setminus K$ carries an $F$-structure 
outside a bounded number of $\varepsilon_{j}$-balls, with $\varepsilon_{j} \rightarrow 
0$ as $\varepsilon \rightarrow 0$, for $K \subset \Omega$ sufficiently large. 

  It is clear that the metric $g_{0}$ on $V$ or $\Omega$ satisfies the Euler-Lagrange equations
$$\nabla{\mathcal F}  = 0 $$
given by (4.19) or (4.33). Further, it is also clear from the construction that the 
configurations $(\Omega_{\varepsilon}, g_{\varepsilon})$ are pointed Gromov-Hausdorff 
limits of sequences of unit volume metrics on $M$. 
{\endproof}

  Next consider the functionals ${\mathcal W}^{2}$, or ${\mathcal W}_{\pm}^{2}$. Although 
we conjecture that Theorem 4.10 also holds for these functionals, there is a basic 
obstacle to proving this. Namely, given an (arbitrary) compact 4-manifold, it is 
unknown if there exists a minimizing sequence $\{g_{i}\}$ for ${\mathcal W}^{2}$ such that 
\begin{equation} \label{e4.54}
\int_{M}s_{g_{i}}^{2} \leq \Lambda,
\end{equation}
for some (arbitrary) $\Lambda < \infty$. Since ${\mathcal W}^{2}$ is conformally invariant, 
one may choose $g_{i}$ to be Yamabe metrics, in which case (4.54) is equivalent to 
\begin{equation} \label{e4.55}
s_{g_{i}} \geq -\Lambda' > - \infty.
\end{equation}
If $M$ has no such minimizing sequence, there seems little hope (at present) of proving 
the existence of a generalized metric realizing $\inf {\mathcal W}^{2}$. On the other hand, 
we conjecture that (4.54) always holds. Of course the bound (4.54) or (4.55) is equivalent 
to a bound on ${\mathcal R}^{2}$, given a bound on ${\mathcal W}^{2}$, via the Chern-Gauss-Bonnet 
theorem. 

  Thus, in the following, we essentially assume (4.54). More precisely, since we are working 
with special minimizing sequences $(\Omega_{\varepsilon,\lambda}, g_{\varepsilon,\lambda})$ 
obtained by minimizing ${\mathcal W}_{\varepsilon,\lambda}^{2p}$, we impose (4.54) on 
$(\Omega_{\varepsilon,\lambda}, g_{\varepsilon,\lambda})$, for a suitable choice of 
$\varepsilon$, $\lambda$. Due to the conformal invariance of the $\varepsilon = 0$ 
limit, one should impose a gauge condition which minimizes ${\mathcal S}^{2}$ in its conformal 
class, (see again Example 3.1). As will seen below, this is implied by the condition that 
as $i \rightarrow \infty$, 
\begin{equation} \label{e4.56}
\varepsilon_{i} \rightarrow 0, \ \  \lambda_{i} \rightarrow 0 \ \ 
{\rm and} \ \  \varepsilon_{i} / \lambda_{i} \rightarrow 0.
\end{equation}
One may either keep $p$ fixed, or let $p_{i} \rightarrow 1$. A resulting sequence 
$(\Omega_{i}, g_{i}) = (\Omega_{\varepsilon_{i},\lambda_{i}}, g_{\varepsilon_{i},\lambda_{i}})$ 
is then called {\it preferred} if there exists a constant $\Lambda < \infty$ such that 
\begin{equation} \label{e4.57}
\int_{\Omega_{i}}s_{g_{i}}^{2} \leq \Lambda .
\end{equation}
The existence of a preferred minimizing sequence clearly depends only on the diffeomorphism 
type of $M$. In fact, it is easily seen to be equivalent to the following condition: 
for each $i$, there exists $\delta_{i}$, with $\delta_{i} \rightarrow 0$ as $i \rightarrow 
\infty$, and a metric $g_{i}$ on $M$ such that (4.54) holds and 
\begin{equation} \label{e4.58}
{\mathcal F}(g_{i}) \leq \inf_{{\mathbb M}_{1}}{\mathcal F}_{\varepsilon_{i},\lambda_{i}} + 
\delta_{i},
\end{equation}
where ${\mathcal F} = {\mathcal W}^{2}$ or ${\mathcal W}_{\pm}^{2}$. 
\begin{theorem} \label{t4.11}
Let $M$ be as in Theorem 4.10 and let ${\mathcal F}$ be one of the functionals ${\mathcal W}^{2}$, 
${\mathcal W}_{\pm}^{2}$. Then any preferred sequence $(\Omega_{i}, g_{i})$ has a subsequence 
converging in the pointed Gromov-Hausdorff topology to a minimizing configuration 
$(\Omega, g_{0})$ for ${\mathcal F}$ which satisfies all the properties of Theorem 4.10. 
Within the conformal class $[g_{0}]$, the metric $g_{0}$ minimizes ${\mathcal S}^{2}$ 
among compact perturbations. The limit $(V, g_{0})$ or $(\Omega, g_{0})$ satisfies 
the Bach equations (4.41).
\end{theorem}

{\bf Proof:} Since the sequence $(\Omega_{i}, g_{i})$ is preferred, it has a uniform 
bound on the $L^{2}$ norm of curvature $R$, cf. (4.57)-(4.58) above. The proof that a 
subsequence converges to a limit $(V, g_{0})$ or $(\Omega, g_{0})$ satisfying the 
properties in Theorem 4.10 is then exactly the same as the proof of Theorem 4.10. 

  The conformal gauge condition satisfied by $g_{0}$ is determined by the form of the 
perturbation ${\mathcal W}_{\varepsilon,\lambda}^{2p}$; different conformal gauges are 
obtained by choosing different perturbations. To obtain the conformal gauge condition for 
$g_{0}$, divide the trace equation (4.40) by $\varepsilon$. A small computation then gives
\begin{eqnarray} \label{e4.59}
-\frac{\lambda}{\varepsilon}\Delta s - 2p\Delta \{(1 + |Ric|^{2})^{p-1}s\} - 
p\langle E, D^{2}(1 + |Ric|^{2})^{p-1} \rangle \\
+ 2(1 + |Ric|^{2})^{p-1}[1 + (1-p)|Ric|^{2} + \frac{c_{\varepsilon}}{\varepsilon}] = 0. 
\nonumber
\end{eqnarray}
It then follows from (4.56) and the smooth convergence of $g_{\varepsilon}$ to $g_{0}$ on 
the regular part $\Omega_{0}$ of $\Omega$, that the metric $g_{0}$ satisfies the limit 
equation 
\begin{equation} \label{e4.60}
\Delta s = 0.
\end{equation}

  In fact, we claim that $g_{0}$ minimizes the functional ${\mathcal S}^{2}$ 
in the conformal class $[g_{0}]$ among comparison metrics which agree with 
$g_{0}$ outside some compact set. To see this, fix any $\varepsilon > 0$, 
$\lambda > 0$, (and $p > 1$), and let $g_{k} = g_{k}(\varepsilon, \lambda, p)$ be 
a sequence of metrics on $M$ minimizing ${\mathcal W}_{\varepsilon,\lambda}^{2p}$ and 
converging to $(\Omega_{\varepsilon,\lambda},g_{\varepsilon, \lambda})$. Now let 
$[g_{k}]$ be the conformal class of $[g_{k}]$ and let $\gamma_{k}$ be metrics in 
$[g_{k}]$ minimizing ${\mathcal W}_{\varepsilon,\lambda}^{2p}$ in the conformal class 
$[g_{k}]$. (If $\gamma_{k}$ is not unique, one may choose a $\gamma_{k}$ 
closest to $g_{k}$ in a given smooth topology). By a compactness result 
of Gursky [21], for each $k$, the metrics $\gamma_{k}$ exist on $M$, and are 
in $[g_{k}]$. Thus, the metrics $\gamma_{k}$ minimize the functional 
$$\frac{2}{\lambda}{\mathcal W}^{2} + \frac{\varepsilon}{\lambda}\int(1 + |Ric|^{2})^{p} 
+ \int s^{2},$$
in $[g_{k}] \subset {\mathbb M}_{1}$. Since ${\mathcal W}^{2}$ is conformally invariant, 
$\gamma_{k}$ minimizes the functional $\frac{\varepsilon}{\lambda}\int(1 + |Ric|^{2})^{p} 
+ \int s^{2}$ in the conformal class $[g_{k}] \subset {\mathbb M}_{1}$.

  Now choose a sequence $\varepsilon_{j} \rightarrow 0$, $\lambda_{j} \rightarrow 0$ 
satisfying (4.56), with $g_{\varepsilon_{j}} \rightarrow  g_{0}$. Then a suitable 
diagonal subsequence of the double sequence $(k, j)$ gives the convergence of 
$\gamma_{k_{j}}$ to $g_{0}$, which proves the claim. Of course, (4.60) is the Euler-Lagrange 
equation for ${\mathcal S}^{2}$ in a given conformal class. This completes the proof of 
Theorem 4.11. 
{\endproof}

\begin{remark} \label{r 4.12.}
  {\rm The following remarks pertain to both Theorems 4.10 and 4.11. It is not 
known if the metric $g_{0}$ is complete on $\Omega$, or if $\Omega$ has 
finitely many components. As in Conjecture 4.7, we conjecture both of these 
to be the case, and that equality holds in (4.44) as well as in (4.45) when 
the analogous contributions of the singularities in the collapsing region 
of $M \setminus K$ are added to the left hand side. Moreover, we conjecture 
that $\Omega$ has finite topological type, and the regular set 
$\Omega_{0}$ embeds in $M$, $\Omega_{0} \subset M$. 

  We point out that it follows from Theorem 4.18 below that all but finitely 
many components of $\Omega$ have an $F$-structure. 

 It is clear that if the curvature $|R|$ of $g_{0}$ is pointwise 
bounded, then $(\Omega , g_{0})$ is complete. In fact, as in Remark 4.3(i), 
completeness follows from just a uniform lower bound on the Ricci curvature, 
$Ric_{g_{0}} \geq  -\lambda$, for some $\lambda < \infty$, by the 
Bishop-Gromov volume comparison theorem. 

 Again, the functionals ${\mathcal W}^{2}$ and ${\mathcal W}_{\pm}^{2}$ are 
conformally invariant, and so completeness is understood to be with 
respect to a metric $g_{0}$ minimizing ${\mathcal S}^{2}$ in its 
conformal class. }
\end{remark}

 Theorems 4.10 and 4.11 give an existence result for idealized minimizers of 
${\mathcal F}$ on a given 4-manifold $M$. As is to be expected, it does not 
assert the existence of metrics on the given (arbitrary) 4-manifold $M$ which 
minimize ${\mathcal F}$; instead the objects are ``generalized metrics'', as 
described by Cases (I) and (II). Thus, these results define a generalized 
moduli space $\widetilde{\mathcal P}_{{\mathcal F}}$ of minimizers, 
associated to each functional ${\mathcal F}$ in (3.1) and 4-manifold $M$.

 By their construction in Theorem 4.10, the moduli spaces 
$\widetilde{\mathcal P}_{{\mathcal F}}$ are compact, in that given any sequence 
$\{g_{i}\} \in  \widetilde{\mathcal P}_{{\mathcal F}}$, a subsequence converges in the 
pointed Gromov-Hausdorff topology to a limit $g \in  
\widetilde {\mathcal P}_{{\mathcal F}}$. The same result holds for the conformally invariant 
functionals in Theorem 4.11, for the portion $\widetilde{\mathcal P}_{{\mathcal F}}(\Lambda)$ 
of the moduli space $\widetilde{\mathcal P}_{{\mathcal F}}$ for which 
$$\int_{\Omega}s^{2} \leq \Lambda,$$
for any given $\Lambda < \infty$. (As seen in Example 3.1, 
$\widetilde{\mathcal P}_{{\mathcal F}}(\Lambda)$ may be a strict subset of 
$\widetilde{\mathcal P}_{{\mathcal F}}$, for all $\Lambda$). 

\medskip

  On the other hand, let ${\mathcal M}_{{\mathcal F}}$ be the ``actual'' 
moduli space of minimizers of ${\mathcal F}$ on $M$, i.e. the space of smooth 
metrics on $M$, modulo diffeomorphisms, which minimize ${\mathcal F}$. By the 
elliptic estimates used in the proof of Proposition 4.6, it is straightforward 
to see that any weak $L^{2,q}$ solution of the Euler-Lagrange equations 
$\nabla {\mathcal F} = 0$, with $q > 4$, is $C^{\infty}$ smooth, (away from 
the orbifold singularities). 

  Of course ${\mathcal M}_{{\mathcal F}}$ may well often be empty. At best, 
the space $\widetilde{\mathcal P}_{{\mathcal F}}$ can only be considered empty 
if $\Omega  = \emptyset$, for all minimizers $(\Omega, g_{0})$ of ${\mathcal F}$, 
or if $M$ has no preferred minimizing sequences in the case of the conformally 
invariant functionals. Even in the case $\Omega = \emptyset$, the space 
${\widetilde{\mathcal P}}_{\mathcal F}$ is better thought of as consisting of 
fully degenerate metrics on $M$ instead of empty; one expects that the fact that 
$\Omega = \emptyset$ leads to strong topological restrictions on $M$. 

  However, perhaps surprisingly, it is not clear that one always necessarily has 
\begin{equation} \label{e4.61}
{\mathcal M}_{{\mathcal F}} \subset \widetilde{\mathcal P}_{{\mathcal F}}.
\end{equation}
Namely, the configurations $(\Omega, g_{0}) \in \widetilde{\mathcal P}_{{\mathcal F}}$ 
are constructed by a very specific process of first passing to minimizers for a 
perturbed functional, and then examining their behavior as the perturbation 
parameter ${\varepsilon}$ is taken to 0. From this construction then, it is not clear 
that one has obtained the full moduli space of minimizers of ${\mathcal F}$. 
In particular, if $g$ is a smooth metric on $M$ which minimizes ${\mathcal F}$, 
then one has to consider the issue whether $g$ may be approximated arbitarily 
closely by minimizers $(\Omega_{\varepsilon}, g_{\varepsilon})$ of the 
$\varepsilon$-perturbed functional as $\varepsilon \rightarrow 0$. We 
conjecture that this is always the case, so that in fact (4.61) does hold. 

  In any case, suppose then that 
$${\mathcal M}_{{\mathcal F}} \neq \emptyset,$$
so that (4.61) is a non-trivial statement, and let $g$ be any smooth metric on $M$ in 
${\mathcal M}_{{\mathcal F}}$. Then it is obvious that $g$ can be approximated by smooth 
metrics on $M$ with a uniform bound on ${\mathcal R}ic^{2p}$. In this context, Theorem 4.1 
and Proposition 4.6 hold, and imply that the completion of ${\mathcal M}_{{\mathcal F}}$ 
with respect to this norm consists of ${\mathcal M}_{{\mathcal F}}$ together with limits 
which have exactly the same structure as described in Cases (I) and (II) of Theorems 4.10 
and 4.11. 

  The following Lemma is useful and of some interest in this setting. 
\begin{lemma} \label{l4.13}
Given $M$ and ${\mathcal F} = {\mathcal R}ic^{2}$ or ${\mathcal R}^{2}$, suppose that 
\begin{equation} \label{e4.62}
{\mathcal M}_{\mathcal F} \neq \emptyset.
\end{equation}
Then minimizers $(\Omega_{\varepsilon}, g_{\varepsilon})$ of the perturbed functional 
${\mathcal F}_{\varepsilon}$ satisfy
\begin{equation} \label{e4.63}
{\mathcal R}ic^{2p}(\Omega_{\varepsilon}, g_{\varepsilon})  \leq \Lambda ,
\end{equation}
for some fixed $\Lambda < \infty$, independent of $\varepsilon$ or the choice of minimizer 
$(\Omega_{\varepsilon}, g_{\varepsilon})$. 
\end{lemma}

{\bf Proof:}
We give the proof in the case of ${\mathcal F} = {\mathcal R}ic^{2}$; the proof in the 
case of  ${\mathcal R}^{2}$ is the same. Let $\gamma \in {\mathcal M}_{\mathcal F}$ be a smooth 
metric on $M$ minimizing ${\mathcal F}$ and let $(\Omega_{\varepsilon}, g_{\varepsilon})$ 
be any minimizing configuration for ${\mathcal F}_{\varepsilon} = {\mathcal R}ic_{\varepsilon}^{2p}$. 
Then $(\Omega_{\varepsilon}, g_{\varepsilon})$ also minimizes the functional 
$$I(g) = \frac{1}{\varepsilon}\int |Ric|^{2} + \int (1 + |Ric|^{2})^{p} - 
\inf_{{\mathbb M}_{1}}\frac{1}{\varepsilon}\int |Ric|^{2},$$
which differs from the functional $\frac{1}{\varepsilon}{\mathcal R}ic^{2p}$ by a constant. 
Let $g_{j}$ be a sequence of metrics on $M$ with ${\mathcal R}ic^{2p}(g_{j})$ uniformly 
bounded, and converging to $(\Omega_{\varepsilon}, g_{\varepsilon})$ in the pointed 
Gromov-Hausdorff topology. Since the metric $\gamma$ realizes 
$\inf_{{\mathbb M}_{1}}\frac{1}{\varepsilon}\int |Ric|^{2}$, one has for any $j$ sufficiently 
large, 
$$I(g_{j}) \leq I(\gamma) + \delta_{j} = \int (1 + |Ric|^{2}(\gamma))^{p} + \delta_{j},$$
where $\delta_{j} \rightarrow 0$ as $j \rightarrow \infty$. The result then follows from 
the obvious inequality $\int (1 + |Ric|^{2}(g_{j}))^{p} \leq I(g_{j})$, and taking the limit 
$j \rightarrow \infty$. 
{\endproof}

  Thus, the hypothesis (4.62) implies that all elements in $\widetilde{\mathcal P}_{\mathcal F}$, 
for ${\mathcal F} = {\mathcal R}ic^{2}$ or ${\mathcal R}^{2}$, are also limits of sequences of 
smooth metrics on $M$ which have a uniform bound on ${\mathcal R}ic^{2p}$. This leads 
naturally to the following definition. 

\begin{definition} \label{d4.14}
{\rm For ${\mathcal F}$ any of the functionals in (3.1), the generalized moduli space 
$\widetilde{\mathcal M}_{\mathcal F}$ is defined to be the space of minimizing 
configurations of ${\mathcal F}$, satisfying the conclusions of Theorems 4.10 or 4.11, 
endowed with the pointed Gromov-Hausdorff topology. }
\end{definition}

  Thus, by definition, one has
\begin{equation} \label{e4.64}
{\mathcal M}_{\mathcal F} \subset \widetilde{\mathcal M}_{\mathcal F}.
\end{equation}
Further, Theorem 4.10 implies that for either of the functionals ${\mathcal R}ic^{2}$ or 
${\mathcal R}^{2}$, the generalized moduli space $\widetilde{\mathcal M}_{\mathcal F}$ 
always exists, while for ${\mathcal W}^{2}$ or ${\mathcal W}_{\pm}^{2}$, Theorem 4.11 
implies the space $\widetilde{\mathcal M}_{\mathcal F}$ exists provided (4.57) or (4.58) 
holds for the manifold $M$. 

 This discussion leads now quite easily to the following structure theorem on the moduli 
space ${\mathcal M}_{\mathcal F}$. This result generalizes the results of [8], [45] 
discussed previously in \S 3(II), cf. also Remark 4.17 below. 

\begin{theorem} \label{t 4.15.}
  Let ${\mathcal F} = {\mathcal R}ic^{2}$ or ${\mathcal R}^{2}$. Then the completion 
$\overline {\mathcal M}_{{\mathcal F}}$ of ${\mathcal M}_{{\mathcal F}}$ with respect to the 
pointed Gromov-Hausdorff topology is contained in $\widetilde {\mathcal M}_{{\mathcal F}}$, 
so that the metrics are of the form described by (I), (II) in Theorem 
4.10. Moreover, the completion $\overline{\mathcal M}_{{\mathcal F}}$ is compact. 

  The same result holds for ${\mathcal F} = {\mathcal W}^{2}$ or ${\mathcal W}_{\pm}^{2}$ for 
the portion ${\mathcal M}_{\mathcal F}(\Lambda)$ given by (4.54), for any given $\Lambda < \infty$. 
\end{theorem}
{\bf Proof:} 
For the cases ${\mathcal F} = {\mathcal R}ic^{2}$ or ${\mathcal R}^{2}$, Lemma 4.13 gives a uniform 
bound on ${\mathcal R}ic^{2p}$ on $\widetilde {\mathcal M}_{{\mathcal F}}$, and so the result 
follows from Theorem 4.1. 

  For the case ${\mathcal F} = {\mathcal W}^{2}$ or ${\mathcal W}_{\pm}^{2}$, Lemma 4.13 does not 
hold. However, one can apply, essentially word for word, the proof of Theorem 4.10 
to sequences $g_{i} \in {\mathcal M}_{\mathcal F}(\Lambda)$, using the elliptic regularity 
estimates from Proposition 4.6, and the result follows from this proof. 
{\endproof}

  There are course many questions one can now begin to ask regarding the 
structure of metrics in $\widetilde{\mathcal M}_{{\mathcal F}}$ and the 
structure of $\widetilde{\mathcal M}_{{\mathcal F}}$ as a whole. Obviously, the 
main question is how the geometric decomposition with respect to ${\mathcal F}$ 
is related to the topology of $M$, as in the Thurston decomposition 
in dimension 3 discussed in \S 2. We will discuss the case of orbifold 
and cusp singularities separately.

\medskip

 As in the Einstein case discussed in \S 3 the geometry and topology 
of the orbifold singularities is quite well understood. Thus, with each 
orbifold singularity $q \in V$ or $q \in \Omega$, one has associated a 
finite number of complete ALE spaces $(N_{j}, \gamma_{j})$ which arise 
as blow-up limits at different scales. Each space $(N_{j}, \gamma_{j})$ 
is a solution of the Euler-Lagrange equations $\nabla{\mathcal F} = 0$, (in 
fact a minimizer for ${\mathcal F}$ among compact perturbations). The 
orbifold singularities at $q$ arise by crushing all of these spaces to 
points, at different scales, as $\varepsilon \rightarrow 0$. 

  As a simple example, consider the Schwarzschild metric $g_{S}$ on 
$S^{3}\times {\mathbb R} \sim  {\mathbb R}^{4}\setminus \{0\}$, given by 
$$g_{S} = (1+\frac{2m}{r^{2}})^{2}g_{Eucl}, $$
where $g_{Eucl}$ is the Euclidean metric on ${\mathbb R}^{4}\setminus \{0\}$. 
This metric is conformally flat, with zero scalar curvature, and so is a 
critical point (in fact a minimizer) of all the functionals ${\mathcal F}$ in 
(3.1). The metric $g_{S}$ has two ALE (in fact AE) ends and the blow-down of 
$g_{S}$ is the orbifold consisting of two copies of ${\mathbb R}^{4}$ glued together 
at the origin $\{0\}$. The same process can be carried out for instance 
for conformally flat metrics with any finite number of ALE ends. Similarly, 
one may have a pair of Schwarzschild metrics joined at a point, resulting 
from blowing-down a Schwarzschild neck joining them at a higher scale. 

 This picture holds in general, and the orbifold singularities 
correspond to a generalized connected sum decomposition of $M$. In 
more detail, suppose one is in Case (I) of Theorem 4.10/4.11. Then 
$M$ is the union of the regular set $V_{0}$ with a finite collection of 
ALE spaces:
$$M = V_{0}\cup\{N_{k}(m)\}, $$
where $V_{0}$ has a finite number of connected components, and a finite 
number of singular points $q_{k}$, $1 \leq k \leq Q$. The union takes place 
along spherical space forms $S^{3}/\Gamma$. We recall, as in 
the example above that the singular points may be topologically regular, 
i.e.~the local fundamental groups attached to $q_{k}$ may be trivial. 
For each fixed singular point $q_{k}$, one has a finite tree of orbifold 
ALE spaces, which are smooth manifolds at the top level. Thus, with any 
$q_{k}$ itself is associated a finite collection of ALE orbifolds 
$N_{k}^{j_{1}}(1)$, whose blow-downs give rise to the singularity at $q_{k}$. 
Each orbifold singularity $q_{k,\ell}^{j_{1}}(1)$ in any $N_{k}^{j_{1}}(1)$ arises 
in the same manner, so that a next level collection of ALE orbifolds 
$N_{k,\ell}^{j_{1}, j_{2}}(2)$ is associated with each $q_{k,\ell}^{j_{1}}(1)$ 
by blowing down. This process repeats itself a finite number of times, until 
reaching the top scale, where the ALE spaces are all smooth manifolds; 
cf. [8], [9] and [11] for further details. 

 A minimizing sequence $g_{i}$ for ${\mathcal F}$ converging to $(V, g_{0})$ 
collapses the full collection of ALE spaces $N_{k}(m)$ to the singular 
points of $V$, (at varying scales). Thus, one may write, (somewhat loosely), 
\begin{equation} \label{e4.65}
M = V\#_{\Gamma}\{N_{k}(m)\},
\end{equation}
where the connected sum is along spherical space-forms, (not just along 3-spheres). 
In most cases, the decomposition (4.65) is topologically non-trivial. However, its 
an open question whether some component in (4.65) could be a 4-sphere; (especially for 
minimizers of ${\mathcal F}$, this seems unlikely). The same description as above holds 
for the orbifold singularities in $\Omega$ in the situation of Case (II) of 
Theorem 4.10/4.11. 

\medskip

On the other hand, essentially nothing is currently known about the relation of the 
topology of $\Omega$ to $M$ when $\Omega \neq V$. In dimension 3, $\Omega$ corresponds 
to the hyperbolic piece $H$, and so is topologically essential in $M^{3}$, 
in that $\pi_{1}(H)$ injects in $\pi_{1}(M^{3})$. It would be very interesting 
to see any progress on this question. 

\medskip

   Continuing with the discussion on orbifolds, let 
$\widetilde M_{{\mathcal F}}^{~o}$ be the part of the generalized moduli 
space $\widetilde {\mathcal M}_{{\mathcal F}}$ consisting of the orbifold 
metrics $(V, g_{0})$, i.e. Case (I) of Theorem 4.10/4.11. Suppose 
$$\widetilde {\mathcal M}_{{\mathcal F}}^{~o} \neq  \emptyset . $$
It is natural to ask if $\widetilde M_{{\mathcal F}}^{~o}$ arises as the 
orbifold frontier $\partial_{o}{\mathcal M}_{{\mathcal F}}$ in the completion 
$\overline{\mathcal M}_{{\mathcal F}}$ of the moduli space ${\mathcal M}_{{\mathcal F}}$ 
of smooth minimizers $(M, g_{0})$ on $M$. In other words, does there 
exist a sequence of smooth metrics $g_{i} \in {\mathcal M}_{{\mathcal F}}$ on $M$, 
for which $(M, g_{i}) \rightarrow (V, g_{0})$ in the 
Gromov-Hausdorff (or $L^{2})$ topology? Thus, the metrics $g_{i}$ on $M$
resolve the orbifold singularities on $(V, g_{0})$. This of course implies 
in particular that 
$${\mathcal M}_{{\mathcal F}} \neq  \emptyset  . $$

 This issue is essentially equivalent to the problem of glueing 
together smooth or orbifold singular metrics in 
$\widetilde {\mathcal M}_{{\mathcal F}}^{~o}$. Namely, one would like to reverse 
the generalized connected sum decomposition (4.65) by a glueing process, 
as follows. Suppose $M$ is a manifold of the form
$$M = V_{0}\cup_{\Gamma}\{N_{k}\}, $$
where the regular set $V_{0}$ of $V$ may have several, although finitely many 
components  and $\{N_{k}\}$ is a collection of ALE manifolds. Suppose $V$ carries 
an orbifold metric $g_{0}$ in $\widetilde{{\mathcal M}_{\mathcal F}}$ and each $N_{k}$ 
carries a complete smooth metric $g_{k}$ minimizing ${\mathcal F}$ on compact sets. 
Can one find a smooth metric $g$ on $M$, minimizing ${\mathcal F}$, which is close 
to the metric $g_{0}$ on $V_{0}$ and to blow-down rescalings of each $g_{k}$ 
on $N_{k}$? 

\medskip

 While this has been a long-standing open question for general Einstein 
metrics on 4-manifolds, quite a lot is now understood on this question 
in the case of self-dual metrics. As mentioned in the beginning of \S 4, 
Taubes [40] proved the first such geometric glueing theorem in the context 
of self-dual solutions of the Yang-Mills equations. General glueing techniques 
and results for self-dual metrics were developed by Floer [20], Donaldson-Friedman [19], 
culminating in the general result of Taubes [41] that given any compact 
oriented 4-manifold $M$, the manifold $\hat M = M\#n{\mathbb C}{\mathbb P}^{2}$ 
for $n$ sufficiently large, (depending on $M$) admits a smooth self-dual metric. 
Explicit self-dual metrics on $n{\mathbb C}{\mathbb P}^{2}$ were first constructed by 
LeBrun [27] and later by Joyce [24]. A very interesting and quite general glueing result 
for self-dual orbifolds has been obtained by Kovalev-Singer [26]; this answers 
the question above in the self-dual case affirmatively in many cases.

\begin{remark} \label{r 4.16.}
  {\rm It is of interest to understand whether the results above hold for 
the much weaker functionals ${\mathcal S}^{2}$ or $-{\mathcal S}|_{{\mathcal Y}}$. 
Thus, consider for instance the perturbation
\begin{equation} \label{e4.66}
{\mathcal S}^{2p}_{\varepsilon} = \int s^{2} + \varepsilon\int 
(1+|Ric|^{2})^{p}. 
\end{equation}
Theorem 4.1 holds of course for ${\mathcal S}^{2p}_{\varepsilon}$.

 The Euler-Lagrange equations of ${\mathcal S}^{2p}_{\varepsilon}$ on 
${\mathbb M}_{1}$ are:
\begin{equation} \label{e4.67}
D^{*}Dh - 2\delta^{*}\delta h -\delta\delta h g -4\delta^{*}\delta Ric 
- 4(\delta\delta Ric) g + {\mathcal P}_{S} = 0, 
\end{equation}
where $h = \varepsilon p(1+|Ric|^{2})^{p-1}Ric$ and ${\mathcal P}_{S}$ is a 
curvature term, given by
$${\mathcal P}_{S} = -2sRic + \tfrac{1}{2}s^{2}g - 2R(h) + 
\tfrac{1}{2}[\varepsilon (1+|Ric|^{2})^{p}+ c]g. $$
Using the same methods as above, it is not difficult to show, although 
we will not carry out the details, that Proposition 4.6 holds for 
${\mathcal S}^{2p}_{\varepsilon}$, i.e. weak solutions of the Euler-Lagrange 
equations are $C^{\infty}$ smooth.

 However the proof of Theorem 4.10 does not hold for ${\mathcal S}^{2}$ (or 
$-{\mathcal S}|_{{\mathcal Y}}$), and one does not expect this result to hold 
for such weak functionals. First, a bound on ${\mathcal S}^{2}$ does not 
imply a bound on ${\mathcal R}^{2}$, and so one immediately loses the 
statements on finiteness of the number of singular points in Theorem 4.10. 
Moreover, the dominant terms in the Euler-Lagrange equation (4.67), i.e. 
$D^{*}Dh$ and $\delta^{*}\delta h$ tend to 0 weakly as 
$\varepsilon \rightarrow 0$, and so the argument in the proof of Theorem 4.10 
concerning smooth convergence to the limit no longer holds. (The equation 
(4.67) is no longer uniformly elliptic as $\varepsilon \rightarrow 0$). Using 
the methods in [7], it may be possible to obtain such smooth convergence in the 
regions where $s \leq  -\lambda$, for some $\lambda > 0$, but this is likely to 
fail in regions where $s \geq 0$. }
\end{remark}

\begin{remark} \label{r4.17}
{\rm We point out that most of the results above do not actually require the limit 
metrics to be minimizers of ${\mathcal F}$. For example, with the exception of the 
statements regarding the infimum in (4.43) and (4.45), Theorem 4.10 holds for the 
components of the moduli space of critical points which are local minimizers of 
${\mathcal F}$ for which 
\begin{equation} \label{e4.68}
{\mathcal F} \leq \Lambda,
\end{equation}
for some $\Lambda < \infty$. (Recall that ${\mathcal F}$ is constant on connected components 
of $\widetilde{\mathcal M}_{\mathcal F}$). Similarly, Theorem 4.11 holds for the part 
$\widetilde{\mathcal M}_{\mathcal F}(\Lambda)$ of $\widetilde{\mathcal M}_{\mathcal F}$ satisfying 
(4.68) and (4.54). 

  For example, all of Theorem 4.15 holds for these parts of the moduli space of 
critical points of ${\mathcal F}$. We also conjecture that the restriction that the 
critical metrics locally minimize ${\mathcal F}$ is not necessary. }
\end{remark}

 We conclude with a simple application of the methods discussed above. 
A well-known question of Gromov asks if there is an $\varepsilon_{0} = 
\varepsilon_{0}(n) > 0$ such that $M$ is a closed $n$-manifold admitting a 
metric such that 
\begin{equation} \label{e4.69}
{\mathcal R}^{n/2} = \int_{M}|R|^{n/2} \leq  \varepsilon_{0}, 
\end{equation}
then is $\inf {\mathcal R}^{n/2} = 0$? In other words, is there a gap for 
the values of $\inf {\mathcal R}^{n/2}$ about 0. The $L^{\infty}$ version of 
this question, i.e. 
\begin{equation} \label{e4.70}
\inf_{\mathbb M} (vol^{2/n}|R|_{L^{\infty}}) \leq  \varepsilon_{0} \Rightarrow  
\inf_{\mathbb M} (vol^{2/n}|R|_{L^{\infty}}) = 0, 
\end{equation}
was proved in dimension 4 by Rong [37], by showing that the hypothesis 
in (4.70) implies the existence of a polarized $F$-structure; the 
conclusion in (4.70) then follows from the work of Cheeger-Gromov [15]. 

 The following result gives a partial answer to Gromov's question in 
dimension 4. 

\begin{theorem} \label{t 4.18. }
  There is an $\varepsilon_{0} > 0$, such that if $M$ is a 4-manifold 
admitting a metric with 
\begin{equation} \label{e4.71}
\int_{M}|R|^{2} \leq  \varepsilon_{0}, 
\end{equation}
then $M$ has an $F$-structure. 
\end{theorem}
{\bf Proof:}
 This is a simple consequence of the proof of Theorem 4.10, with which we thus 
assume some familiarity. Let $(\Omega, g_{0})$ be a minimizing configuration 
for ${\mathcal R}^{2}$, given by Theorem 4.10. Suppose first $\Omega = \emptyset$. 
This may happen in two ways. First the minimizing configurations $(\Omega_{\varepsilon}, 
g_{\varepsilon})$ for ${\mathcal R}_{\varepsilon}^{2p}$ may be non-empty, but collapse 
everywhere as $\varepsilon \rightarrow 0$. Second, $(\Omega_{\varepsilon}, g_{\varepsilon})$ 
may be empty for any $\varepsilon$ sufficiently small. In either case, it follows from 
Theorem 4.10 that $M$ has an $F$-structure metrically on the complement of finitely 
many singular points. The singular points arise from a concentration of curvature in 
$L^{2}$. However, Lemma 4.2 shows that each singular point contributes a definite amount, 
$c_{1}$, to the integral $\int|R|^{2}$. Hence, (4.71) implies there are no singular points, 
and so $M$ itself has an $F$-structure. 

  Next suppose $\Omega \neq \emptyset$. Exactly the same argument rules out any 
orbifold singular points in $(\Omega, g_{0})$, (or $(V, g_{0})$). To prove then 
that $M$ has an F-structure, it suffices to prove that $\Omega$ has an F-structure, 
since the complementary region $M \setminus K$ already has an F-structure for $K$ 
sufficiently large. In turn, the statement that $\Omega$ has an F-structure follows 
from the claim that there exists $\varepsilon_{1} = \varepsilon_{1}(\varepsilon_{0})$, 
such that 
\begin{equation} \label{e4.72}
(\nu^{2}|R|)(x) \leq \varepsilon_{1},
\end{equation}
for all $x \in (\Omega, g_{0})$; here $\nu$ is the volume radius and $|R|$ is the 
pointwise norm of the curvature. The estimate (4.72) is scale-invariant and in 
the scale $\nu = \nu (x) = 1$ requires $|R|(x) \leq \varepsilon_{1}$. Now the 
metric $g_{0}$ satisfies the Euler-Lagrange equations (4.35) and so by the elliptic 
estimates in Proposition 4.6 or Corollary 4.8, one has in this scale,
$$\sup_{B_{x}(1/2)}|R|^{2} \leq c \int_{B_{x}(1)}|R|^{2},$$
for a fixed constant $c < \infty$, independent of $x$ and $g_{0}$. Since the 
bound (4.71) is also scale-invariant, this proves (4.72), which proves the result. 
{\endproof}

 The proof of Theorem 4.18 above shows in fact that if $M$ has a metric satisfying 
$0 <\int_{M}|R|^{2} \leq \varepsilon_{0}$, then $M$ has a (possibly distinct) 
metric which is $\varepsilon_{1}$-volume collapsed, i.e. 
$vol \leq  \varepsilon_{1}$, with locally bounded curvature, i.e. 
$(inj^{2}|R|)(x) \leq  \varepsilon_{1}$, where ${\varepsilon}_{1} = 
\varepsilon_{1}(\varepsilon_{0})$. Thus, the question (4.69) 
follows if Rong's result can be generalized from a global $L^{\infty}$ 
bound on curvature to a local bound. 

 Via the Chern-Gauss-Bonnet theorem, the condition (4.71) can be 
reexpressed as
$$\chi (M) + \frac{1}{8\pi^{2}}\int_{M}|z|^{2} \leq  \varepsilon_{0}, $$
or 
$$-4\chi (M) + \frac{1}{\pi^{2}}\int_{M}|W|^{2} + \frac{1}{24}s^{2} 
\leq  \varepsilon_{0}. $$
Of course, the existence of an $F$-structure implies that $\chi (M) = 0$.

 We point out that the proof of equality in (4.44), (and in (4.45) when the 
contribution of the singularities in the collapsed part is added), is likely to 
require a positive solution of the question (4.69); in fact the two questions are 
probably equivalent.

\bibliographystyle{plain}

\begin{thebibliography}{WWW}
\bibitem [33]{33} G. Perelman, The entropy formula for the Ricci flow and its geometric 
applications, math.DG/0211159.
\end{thebibliography}

\begin{center}
August, 2005
\end{center}

\smallskip
\noindent
\address{Deptartment of Mathematics\\
S.U.N.Y. at Stony Brook\\
Stony Brook, N.Y. 11794-3651}

\noindent
E-mail: anderson@math.sunysb.edu

\end{document}